\documentclass{amsart}

\usepackage{amscd,amssymb}

\def\Q{{\mathbb Q}}
\def\R{{\mathbb R}}
\def\C{{\mathbb C}}
\def\P{{\mathbb P}}
\def\M{{\mathbb M}}
\def\V{{\mathbb V}}
\def\W{{\mathbb W}}
\def\bL{{\mathbb L}}

\def\B{{\mathcal B}}
\def\F{{\mathcal F}}
\def\G{{\mathcal G}}
\def\I{{\mathcal I}}
\def\L{{\mathcal L}}
\def\O{{\mathcal O}}
\def\U{{\mathcal U}}

\def\cV{{\mathcal V}}
\def\cP{{\mathcal P}}

\def\u{{\mathfrak u}}
\def\m{{\mathfrak m}}
\def\g{{\mathfrak g}}
\def\s{{\mathfrak s}}
\def\p{{\mathfrak p}}

\def\bG{\boldsymbol{\G}}
\def\bP{\boldsymbol{\cP}}
\def\bO{\boldsymbol{\O}}

\def\A{\mathbf{A}}
\def\zbar{{\overline{z}}}
\def\Xbar{{\overline{X}}}
\def\Ybar{{\overline{Y}}}
\def\Vbar{\overline{\mathcal{V}}}
\def\Obar{\overline{\O}}
\def\Deltabar{\overline{\Delta}}

\def\Gdr{\G^{\mathrm{DR}}}			
\def\Udr{\U^{\mathrm{DR}}}			

\def\dt{{\bullet}}
\def\Efin{E_{{\rm fin}}}			
\def\Xtilde{\widetilde{X}}
\def\xtilde{\tilde{x}}
\def\ytilde{\tilde{y}}

\def\alphahat{\widehat{\alpha}}
\def\alphatilde{\tilde{\alpha}}

\def\gammatilde{{\tilde{\gamma}}}
\def\rhotilde{\tilde{\rho}}

\def\omegatilde{\widetilde{\omega}}
\def\what{\widehat{w}}
\def\Uhat{\widehat{U}}
\def\Psihat{\widehat{\Psi}}

\def\blank{\phantom{x}}
\def\omit{\underline{\blank}}

\newcommand\spec{\operatorname{Spec}}
\newcommand\im{\operatorname{im}}		

\newcommand\id{\operatorname{id}}

\newcommand\Hom{\operatorname{Hom}}
\newcommand\End{\operatorname{End}}

\newcommand\Aut{\operatorname{Aut}}

\newcommand\Gr{\operatorname{Gr}}
\newcommand\comptensor{\operatorname{\widehat{\otimes}}}

\newtheorem{theorem}{Theorem}[section]
\newtheorem{lemma}[theorem]{Lemma}
\newtheorem{proposition}[theorem]{Proposition}
\newtheorem{corollary}[theorem]{Corollary}

\theoremstyle{definition}
\newtheorem{definition}[theorem]{Definition}
\newtheorem{example}[theorem]{Example}

\theoremstyle{remark}
\newtheorem{remark}[theorem]{Remark}

\begin{document}

\title[The Hodge de~Rham Theory of Relative Completion]%
{The Hodge de~Rham Theory of Relative Malcev Completion}

\author{Richard M.~Hain}

\thanks{This work was partially supported by the National Science
Foundation.}

\address{Department of Mathematics\\ Duke University\\
Durham, NC 27708-0320}

\email{hain@math.duke.edu}

\date{July 18, 1996}

\maketitle

\section{Introduction}

Suppose that $\pi$ is an abstract group, that $S$ is a reductive
algebraic group defined over a field $F$ of characteristic zero,
and that $\rho : \pi \to S(F)$ is a homomorphism with Zariski dense
image. The completion of $\pi$ relative to $\rho$ is a proalgebraic
group $\G$ which is an extension
$$
1 \to \U \to \G \stackrel{p}{\to} S \to 1
$$
where $\U$ is prounipotent, and a homomorphism $\rhotilde : \pi \to \G(F)$
which lifts $\rho$:
$$
\begin{CD}
\pi @>{\rho}>> S \cr
@V{\rhotilde}VV @| \cr
\G @>p>> S
\end{CD}
$$
It is characterized by the following universal mapping property. If
$\phi$ is a homomorphism of $\pi$ to a (pro)algebraic group $G$
over $F$ which is an extension
$$
1 \to U \to G \to S \to 1
$$
of $S$ by a unipotent group $U$, and if the the composite
$$
\pi \to G \to S
$$
is $\rho$, then there is a unique homomorphism $\G \to G$ of
$F$-proalgebraic groups which commutes with the projections to $S$ and
through which $\phi$ factors.

When $S$ is the trivial group, $\G$ is simply the classical Malcev
(or unipotent) completion of $\pi$. In this case, with $F=\R$ or $\C$,
and $\pi$ the fundamental group of a  smooth manifold, there is a de~Rham
theorem for $\O(\G)$ which was proved by K.-T.~Chen \cite{chen}.
In these notes we generalize Chen's de~Rham Theorem from the unipotent
case to the general case. Our approach is based on the notes
\cite{deligne:letter} of Deligne where an approach to computing the
Lie algebra of the prounipotent radical of $\G$ via Sullivan's minimal
models is sketched. Before explaining our result in general, we recall
Chen's de~Rham Theorem in the unipotent case.

If $M$ is a smooth manifold and $w_1,\dots, w_r$ are smooth
1-forms on $M$, then Chen defined
$$
\int_\gamma w_1\dots w_r =
\idotsint\limits_{0\le t_1\le \cdots \le t_r \le 1} f_1(t_1)\dots
f_r(t_r)\, dt_1 \dots dt_r
$$
where $\gamma : [0,1] \to M$ is a piecewise smooth path and
$\gamma^\ast w_j = f_j(t)\, dt$. These are viewed as functions on the
path space of $M$. An iterated integral is a linear combination of
such functions and the constant function. Fix a base point $x\in M$.
Set $\pi = \pi_1(M,x)$. Denote the iterated integrals on the space
of loops in $M$ based at $x$ by $\I_x$. Denote by $H^0(\I_x)$ those
elements of $\I_x$ whose value on a loop depends only on its homotopy
class. Then Chen's $\pi_1$ de~Rham Theorem asserts that integration
induces a Hopf algebra isomorphism
$$
\O(\U) \cong H^0(\I_x)
$$
where $\U$ denotes the real points of the unipotent completion of $\pi$
and $\O(\U)$ its coordinate ring.
Another important ingredient of Chen's theorem is that it gives an
algebraic description of $\I_x$ and $H^0(\I_x)$ in terms of the (reduced)
bar construction on the de~Rham complex of $M$ and the augmentation induced
by the base point.

In this paper we generalize the definition of iterated integrals and
prove a more general de~Rham theorem in which the Hopf algebra $\O(\G)$
of functions on the completion of $\pi_1(M,x)$ relative to a homomorphism
$\rho : \pi_1(M,x) \to S$ is isomorphic to a Hopf algebra of ``locally
constant iterated integrals,'' defined algebraically in terms of a suitable
(2-sided) bar construction on a complex $\Efin^\dt(M,\O(P))$. This complex
of forms plays a central role in all our constructions and was introduced by
Deligne in his notes \cite{deligne:letter}, the main result of which is
that the pronilpotent Lie algebra associated to the 1-minimal model of
$\Efin^\dt(M,\O(P))$ is the Lie algebra of the prounipotent radical
$\U$ of $\G$.

In Section~\ref{groupoid} we define the completion of the fundamental
groupoid of a manifold $M$ with respect to the representation $\rho$.
This is a category (in fact, a groupoid) whose objects are the points
of $X$ and where the Hom sets are proalgebraic varieties; the
automorphism of the object $x\in M$ is the completion of $\pi_1(M,x)$
relative to $\rho$. There is a canonical functor of the fundamental
groupoid of $M$ to this category. We give a de~Rham description of the
coordinate ring of each Hom variety in terms of a suitable 2-sided bar
construction on $\Efin^\dt(M,\O(P))$ and of the functor from the
fundamental groupoid to its relative completion using iterated integrals.

One of the main applications of Chen's $\pi_1$ de~Rham Theorem is to
give a direct construction of Morgan's mixed Hodge structure \cite{morgan}
on the unipotent completion of the fundamental group of a pointed complex
algebraic variety as is explained in \cite{hain:geom}. In this paper we
prove that if $X$ is a smooth complex algebraic variety (or the complement
of a normal crossings divisor in a compact K\"ahler manifold) and
$\V \to X$ is an admissible variation of Hodge structure with polarization
$\langle\blank,\blank\rangle$ whose monodromy representation
$$
\rho : \pi_1(X,x) \to S := \Aut(V_x,\langle\blank,\blank\rangle)
$$
has Zariski dense image\footnote{The assumption that the monodromy
have Zariski dense monodromy can probably be removed. What one needs
to know is that the Zariski closure of the image of $\rho$ is reductive
and that its coordinate ring has a natural real Hodge structure --- see
Remark~\ref{extended}. This should follow from the work of Simpson and Corlette
as each of them has pointed out.}, then the coordinate ring $\O(\G)$ of the
completion of $\pi_1(X,x)$ relative to $\rho$ has a natural mixed
Hodge structure. More generally, we show that the coordinate rings
of the Hom sets of the relative completion of the fundamental groupoid
of $X$ with respect to $\rho$ have canonical mixed Hodge structures.

Our principal application of the Hodge theorem for relative completion
appears in \cite{hain:torelli} where we use it to prove that the unipotent
completion of each Torelli group (genus $\neq 2$) has a canonical mixed
Hodge structure given the choice of a smooth projective curve of genus $g$.
Another application suggested by Ludmil Kartzarkov, and proved in
Section~\ref{hodge_str}, is a generalization
of the theorem of Deligne-Griffiths-Morgan-Sullivan (DGMS) on fundamental
groups of compact K\"ahler manifolds: If $X$ is a compact K\"ahler
manifold and $\V\to X$ is a polarized variation of Hodge structure
with Zariski dense monodromy, then the prounipotent radical of the
completion of $\pi_1(X,x)$ relative to the monodromy representation
has a presentation with only quadratic relations. The theorem of
DGMS is recovered by taking $\V$ to be the trivial variation $\Q_X$.

In Section~\ref{connection} we show that if $X$ is a smooth variety
and $\V$ is an admissible variation of Hodge structure over $X$ with
Zariski dense monodromy representation $\rho$, then there is a canonical
integrable 1-form
$$
\omega \in E^1(X')\comptensor \Gr^W_\dt \u
$$
where $X'$ is the Galois covering of $X$ with Galois group $\im\rho$,
and $\u$ the Lie algebra of the prounipotent radical $\U$ of the
completion $\G$ of $\pi_1(X,x)$ with respect to $\rho$. This
form is $\im \rho$ invariant under the natural actions of $\im \rho$
on $X'$ and $\u$. It can be integrated to the canonical representation
$$
\rhotilde : \pi_1(X,x) \to S\ltimes \U \cong \G.
$$
In the particular case where $X$ is the complement of the discriminant
locus in $\C^n$, where $\pi_1(X,x)$ is the braid group$B_n$ and
$S$ the symmetric group, this connection is the standard one
$$
\omega = \sum_{i<j} d\log(x_i-x_j)\, X_{ij}
$$
on $X'$, the complement in $\C^n$ of the hyperplanes $x_i=x_j$. Kohno
\cite{kohno} used the $\Sigma_n$ invariant form $\omega$ and
finite dimensional representations of $\Gr^W_\dt \u$ to
construct Jones's representations of $B_n$.
Our construction is used in \cite{hain:torelli} to construct an
analogous ``universal projectively flat connection'' for the mapping
class groups in genus $\ge 3$.

I am very grateful to Professor Deligne for sharing his notes on
the de~Rham theory of relative completion with me and for his interest
in this work. I would also like to thank M.~Saito for explaining some
of his work to me, and Hiroaki Nakamura for his careful reading the
manuscript and his many useful comments. I'd also like to thank Kevin
Corlette and Carlos Simpson for freely sharing their ideas on
(\ref{extended}). The bulk of this
paper was written when I was visiting Paris in spring 1995. I would
like to thank the Institute Henri Poincar\'e and the Institute des
Hautes \'Etudes Scientifiques for their generous hospitality and support.

\section{Conventions}

Here, to avoid confusion later on, we make explicit our basic conventions
and review some basic constructions that depend, so some extent, on these
conventions.

Throughout these notes, $X$ will be a connected smooth manifold. By a path
in $X$ from $x\in X$ to $y\in Y$, we shall mean a piecewise smooth map
$\alpha : [0,1] \to X$ with $\alpha(0)=x$ and $\alpha(1) =y$. The set
of all paths in $X$ will be denoted by $PX$. There is a natural
projection $PX \to X\times X$; it takes $\alpha$ to its endpoints
$(\alpha(0),\alpha(1))$.
The fiber of this map over $(x,y)$ will be denoted by $P_{x,y}X$, and the
inverse image of $\{x\}\times X$ will be denoted by $P_{x,-}$. The sets
$PX$, $P_{x,y}X$, $P_{x,-}X$, each endowed with the compact-open
topology, are topological spaces.

We shall multiply paths in their natural order, as distinct from the
functional order. That is, if $\alpha$ and $\beta$ are two paths in $X$
with $\alpha(1) = \beta(0)$, then the path $\alpha\beta$ is defined and
is the path obtained by first traversing $\alpha$, and then $\beta$.

Suppose that $(\Xtilde,\xtilde_o) \to (X,x_o)$ is a pointed universal
covering of $X$. With our path multiplication convention, $\pi_1(X,x_o)$
acts on the {\em left} of $\Xtilde$. One way to  see this clearly is to
note that there is a natural bijection
$$
\coprod_{y\in X}\pi_0(P_{x_o,y}X) \to \Xtilde.
$$
This bijection is constructed by taking the homotopy class of the path
$\alpha$ in $X$ that starts at $x_o$ to the endpoint $\alphatilde(1)$ of
the unique lift $\alphatilde$ of $\alpha$ to $\Xtilde$ that starts at
$\xtilde_o$. With respect to this identification, the action of
$\pi_1(X,x_o)$ is by left multiplication.

Another consequence of our path multiplication convention is that
$\pi_1(X,x_o)$ naturally acts on the {\em right} of the fiber over $x_o$
of a flat bundle over $X$, as can be seen from an elementary computation.
Conversely, if
$$
\rho : F\times \pi_1(X,x_o) \to F
$$
is a right action of $\pi_1(X,x_o)$ on $F$, then one can define
$F\times_\rho \Xtilde$ to be the quotient space $F\times \Xtilde/\sim$,
where the equivalence relation is defined by
$$
(f,gx) \sim (fg,x)
$$
for all $g \in \pi_1(X,x_o)$. This bundle has a natural flat structure
--- namely the one induced by the trivial flat structure on the bundle
$F \times \Xtilde \to \Xtilde$.
The composite
$$
F \cong F \times \{\xtilde_o\} \hookrightarrow
F \times \Xtilde \to F\times_\rho\Xtilde
$$
gives a natural identification of the fiber over $x_o$ with $F$. With
respect to this identification, the monodromy representation of the flat
bundle $F\times_\rho \Xtilde \to X$ is $\rho$.

Of course, left actions can be converted into right actions by using
inverses. Presented with a natural left action of $\pi_1(X,x_o)$ on a
space, we will convert it, in this manner, into a right action in order
to form the associated flat bundle.

The flat bundle over $X$ corresponding to the right $\pi_1(X,x_o)$-module
$V$ will be denoted by $\V$.
For a flat vector bundle $\V$ over $X$, we shall denote the complex of
smooth forms with coefficients in the corresponding $C^\infty$ vector
bundle by $E^\dt(X,\V)$. This is a complex whose cohomology
is naturally isomorphic to $H^\dt(X,\V)$. In particular, the $C^\infty$
de~Rham complex of $X$ will be denoted by $E^\dt(X)$.

By definition, mixed Hodge structures (MHSs) are usually finite
dimensional. When studying MHSs on completions of fundamental
groups, one encounters two kinds of infinite dimensional MHSs
$$
((V_\R,W_\dt),(V_\C,W_\dt,F^\dt)).
$$
In both cases, the weight graded quotients are finite dimensional.
In one, the weight filtration is bounded below (i.e. $W_lV=0$, for some
$l$) so that each $W_mV$ is finite dimensional. In this case we require
that each $W_mV$ with the induced filtrations be a finite dimensional
MHS in the usual sense. The other case is dual. Here the
weight filtration is bounded above (i.e., $V=W_lV$ for some $l$). In
this case, each $V/W_mV$ is finite dimensional. We require that $V$
be complete in the topology defined by the weight filtration (i.e.,
$V$ is the inverse limit of the $V/W_mV$), that each part of the
Hodge filtration be closed in $V$, and that each $V/W_mV$ with the
induced filtrations be a finite dimensional MHS in the usual sense.
Such mixed Hodge structures form an abelian category, as is easily
verified.

Finally, if $V^\dt$ is a graded module and $r$ is a integer, $V[r]^\dt$
denotes the graded module with
$$
V[r]^n = V^{r+n}.
$$

\section{The Coordinate Ring of a Reductive Linear Algebraic Group}
\label{coord}

Suppose that $S$ is a reductive linear
algebraic group over a field $F$ of characteristic zero. The right and
left actions of $S$ on itself induce commuting left and right actions of
$S$ on its coordinate ring $\O(S)$.

If $V$ is a right $S$ module, its dual $V^\ast := \Hom_F(V,F)$ is a left
$S$ module via the action
$$
(s \cdot \phi) (v) := \phi(v\cdot s),
$$
where $s\in S$, $\phi \in \Hom_F(V,F)$ and $v\in V$.

The following result generalizes to reductive groups a well known fact
about the group ring of a finite group.

\begin{proposition}\label{decomp} If $\left(V_\alpha\right)_\alpha$ is a
set of representatives of the isomorphism classes of irreducible
right $S$-modules, then, as an $(S,S)$ bimodule, $\O(S)$ is canonically
isomorphic to
$$
\bigoplus_\alpha V_\alpha^\ast\boxtimes V_\alpha.
$$
\end{proposition}

\begin{proof} This follows from the following facts:
\begin{enumerate}
\item If $V$ is an $S$ module, then the set of matrix entries of
$V$ is the dual $(\End V)^\ast$ of $\End V$. It has commuting right
and left $S$ actions. The right action is induced by left multiplication
of $S$ on itself by left translation, and the left action by the right
action of $S$ on itself.
\item As a vector space, $(\End_F V)^\ast$ is naturally isomorphic to
$V^\ast\otimes V$. The isomorphism takes
$\phi\otimes v\in V^\ast \otimes V$ to the matrix entry
$$
\{f:V \to V\} \mapsto \left\{F \stackrel{v}{\to} V \stackrel{f}{\to}
V \stackrel{\phi}{\to} F \right\}.
$$
It is easily checked that this isomorphism gives an isomorphism
$(\End V)^\ast \cong V^\ast \boxtimes V$ of $(S,S)$-bimodules.
\item By standard arguments (cf.\ \cite{cartier}), the fact that
$S$ is reductive implies that the subspace of $\O(S)$ spanned by
the matrix entries of all irreducible linear representations
is a subalgebra of $\O(S)$. That is, the image of the linear
map
$$
\Phi : \sum_\alpha V_\alpha^\ast \boxtimes V_\alpha \to \O(S)
$$
is a subalgebra of $\O(S)$. Since $\Phi$ is $S\times S$ equivariant, and
since the $V_\alpha^\ast \boxtimes V_\alpha$ are pairwise non-isomorphic
irreducible representations of $S\times S$, $\Phi$ is injective.
\item Since $S$ is linear, it has a faithful linear representation $V_0$,
say and $\O(S)$ is generated by the matrix entries of $V_0$. It follows
that $\Phi$ is surjective, and therefore an algebra isomorphism.
\end{enumerate}
\hfill \end{proof}

Recall that if $G$ is an affine algebraic group, then the Lie
algebra $\g$ of $G$ can be recovered from $\O(G)$ as follows:
Denote the maximal ideal in $\O(G)$ of functions that vanish at
the identity by $\m$. Then, as a vector space, $\g$ is isomorphic to
the dual $\m/\m^2$ of the Zariski tangent space of $G$ at the identity.
The bracket is induced by the comultiplication
$$
\Delta : \O(G) \to \O(G)\otimes \O(G)
$$
as we shall now explain.
Evaluation at the identity and inclusion of scalars give linear
maps $\O(G) \to k$ and $k\to \O(G)$. There is therefore a canonical
isomorphism
$$
\O(G) \cong k \oplus \m.
$$
Using this decomposition, we see that the diagonal induces a
diagonal map
$$
\Deltabar : \m \to \m\otimes \m.
$$
Denote the involution $f\otimes g \mapsto g\otimes f$ of $\m \otimes \m$
by $\tau$. The map
$$
\Deltabar - \tau \circ \Deltabar : \m \to \m\otimes \m
$$
induces the map
$$
\Delta^c : \m/\m^2 \to \m/\m^2\otimes \m/\m^2
$$
dual to the bracket.

\section{A Basic Construction}

{}From this point on $S$ will be a linear algebraic group defined over
$\R$. We will abuse notation and also denote its group of real points by
$S$. We will assume now that we have a representation
$$
\rho : \pi_1(X,x_o) \to S
$$
whose image is Zariski dense.
We will fix a set of representatives $\left(V_\alpha\right)_\alpha$ of
the isomorphism classes of rational representations of $S$.

Composing $\rho$ with the action of $S$ on itself by {\it right}
multiplication, we obtain a right action of $\pi_1(X,x_o)$ on $S$.
Denote the corresponding flat bundle by
$$
p: P \to X.
$$
This is a left principal $S$ bundle whose fiber $p^{-1}(x_o)$ over
$x_o$ comes with an identification with $S$ ; the $S$ action and the
marking of $p^{-1}(x_o)$ are induced by the obvious left action of
$S$ on $S \times \Xtilde$ and by the composite
$$
S \cong S \times \{x_o\} \hookrightarrow S \times \Xtilde \to P.
$$
The point $\xtilde_o$ of $p^{-1}(x_o)$ corresponding to $1 \in S$
will be used as a basepoint of $P$.

Each rational representation of $S$ gives rise to a representation of
$\pi_1(X,x_o)$, and therefore to a local system over $X$. We shall call
such a local system a {\it rational local system}.

The action of $\pi_1(X,x_o)$ on $S$ by right multiplication induces a
left action of $S$ on $\O(S)$, the coordinate ring of $S$. Convert this
to a right action using inverses:
$$
(f\gamma)(s) = f(s\gamma^{-1}),
$$
where $f\in \O(S)$, $\gamma\in \pi_1(X,x_o)$, and $s\in S$. Denote the
associated flat bundle by
$$
\O(P) \to X.
$$
This is naturally a {\em right} flat principal $S$ bundle over $X$. It
follows from (\ref{decomp}) that it is the direct sum of its rational
sub-local systems:
\begin{equation}\label{decomp2}
\O(P) = \bigoplus_\alpha \V_\alpha^\ast \otimes V_\alpha.
\end{equation}
In particular, it is the direct limit of its rational sub-local systems.
Define
$$
\Efin^\dt(X,\O(P)) = \lim_\to E^\dt(X,\M),
$$
where $\M$ ranges over the rational sub-local systems of $\O(P)$.
Denote the cohomology
$$
\lim_\to H^\dt(X,\M)
$$
of this complex by $H^\dt(X,\O(P))$.
The right action of $S$ on $\O(P)$ induces a right action of $S$ on
$$
H^\dt(X,\O(P)).
$$
{}From (\ref{decomp2}), it follows that there is a natural isomorphism
$$
\Efin^\dt(X,\O(P)) \cong
\bigoplus_\alpha E^\dt(X,\V_\alpha^\ast)\otimes V_\alpha
$$
of right $S$ modules. The following result is an immediate consequence.

\begin{proposition}\label{isom}
For each irreducible representation $V$ of $S$, there is a natural
isomorphism
$$
\left[ H^k(X,\O(P))\otimes V\right]^S \cong H^k(X,\V). \qed
$$
\end{proposition}

The bundle $P\to X$ is foliated by its locally flat sections. Denote
this foliation by $\F$. We view it as a sub-bundle of $TP$, the tangent
bundle of $P$. Denote by $E^k(P,\F)$ the vector space consisting
of $C^\infty$ sections of the dual of the bundle
$$
\Lambda^k \F \to P.
$$
One can differentiate sections along the leaves to obtain an
exterior derivative map
$$
d : E^k(P,\F) \to E^{k+1}(P,\F).
$$
With this differential, $E^\dt(P,\F)$ is a differential graded algebra.
Moreover, the left action of $S$ on $P$ induces a natural right action
of $S$ on it, and the natural restriction map
\begin{equation}\label{res}
E^\dt(P) \to E^\dt(P,\F)
\end{equation}
is an $S$-equivariant homomorphism of differential graded algebras.

The base point $\xtilde_o\in P$ induces augmentations
$$
\Efin^\dt(X,\O(P)) \to \R\text{ and } E^\dt(P,\F) \to \R.
$$

\begin{proposition}\label{dga_homom}
There is a natural, augmentation preserving d.g.\ algebra
homomorphism
$$
\Efin^\dt(X,\O(P)) \to E^\dt(P,\F)
$$
which is injective and $S$-equivariant with respect to the natural
right $S$ actions.\hfill \qed
\end{proposition}

\section{Iterated Integrals and Monodromy of Flat Bundles}

Consider the category $\B(X,S)$ whose objects are flat vector
bundles $\V$ over $X$ that admit a finite filtration
$$
\V = \V^0 \supset \V^1 \supset \V^2 \supset \cdots
$$
by sub-local systems with the properties:
\begin{enumerate}
\item the intersection of the $\V^i$ is trivial;
\item each graded quotient $\V^i/\V^{i+1}$ is the local system
associated with a rational representation of $S$.
\end{enumerate}

Denote the fiber over the base point $x_o$ by $V_o$. It has
a filtration corresponding to the filtration $\V^\dt$ of
$\V$:
$$
V_o = V_o^0 \supset V_o^1 \supset V_o^2 \supset \cdots
$$
The second condition above implies that there are rational
representations $\tau_i : S \to \Aut \Gr^i V_o$ such that
the representation of $\pi_1(X,x_o)$ on $\Gr^i V_o$ is the
composite
$$
\pi_1(X,x_o) \stackrel{\rho}{\to} S
\stackrel{\tau_i}{\to} \Aut \Gr^i V_o.
$$
Let $\tau : S \to \prod \Aut \Gr^i V_o$ be the product of the
representations $\tau_i$.
Let
$$
G = \left\{\phi \in \Aut V_o :
\phi\text{ preserves } V_o^\dt\text{ and }
\Gr^\dt\phi\in \im \tau \right\}.
$$
This is a linear algebraic group which is an extension of $\im \tau$
by the unipotent group
$$
U = \left\{\phi \in \Aut V_o : \phi\text{ preserves }
V_o^\dt\text{ and acts trivially on } \Gr^\dt V_o \right\}
$$
whose Lie algebra we shall denote by $\u$.
We shall denote the monodromy representation at $x_o$ of $\V$ by
$$
\rhotilde : \pi_1(X,x_o) \to G.
$$

Denote the $C^\infty$ vector bundles associated to the flat bundles
$\V$ and $\V^i$ by $\cV$ and $\cV^i$, respectively.
We would like to trivialize $\cV$. In order to do this, we pull it
back to $P$ along the projection $p:P \to X$.

\begin{proposition}\label{triv}
There is a trivialization
$$
p^\ast \cV \stackrel{\cong}{\to} P\times V_o
$$
and a splitting of the natural map $G \to \im \tau$ which satisfy
\begin{enumerate}
\item the corresponding connection form%
\footnote{Our convention is that the connection form associated to the
trivialized bundle $V\times X \to X$ with connection $\nabla$ is the
1-form $\omega$ on $X$ with values in $\End V$ which is characterized
by the property that for all sections $f:X \to V$
$$
\nabla f = df - f\omega \in E^1(X)\otimes V.
$$}
$\omegatilde$ satisfies
$$
\omegatilde \in E^1(P)\otimes \u;
$$
\item\label{cond} the isomorphism $V_o \to V_o$, induced by the
trivialization of $p^\ast \cV$ between the fiber over the points $\xtilde_o$
and $s\cdot\xtilde_o$ of $p^{-1}(x_o)$, is $\tau(s)^{-1}$.
\end{enumerate}
\end{proposition}

Note that the second condition implies that the isomorphism $V_o \to V_o$,
induced by the trivialization of $p^\ast \cV$ between the fiber over the
points $a\cdot \xtilde_o$
and $sa\cdot\xtilde_o$ of $p^{-1}(x_o)$, is $\tau(s)^{-1}$.

The first step in the proof is the following elementary result. It
can be proved by induction on the length of the filtration. It gives
the splitting of $G \to \im \tau$.

\begin{lemma}\label{splitting}
There is an isomorphism
$$
\cV \cong \bigoplus_{i\ge 0} \Gr^i \cV
$$
of $C^\infty$ vector bundles that splits the filtration $\cV^\dt$.
That is,
\begin{enumerate}
\item the sub-bundle $\cV^i$ corresponds to $\oplus_{j\ge i} \Gr^j \cV$;
\item the isomorphism
$$
\Gr^i \cV \to \cV^i/\cV^{i+1}
$$
induced by the trivialization is the identity. \hfill \qed
\end{enumerate}
\end{lemma}

\begin{proof}[Proof of (\ref{triv})]
Pulling back the splitting given by (\ref{splitting}) of the
filtration $\cV^i$ to $P$, we obtain a splitting
$$
p^\ast\cV \cong \bigoplus_i p^\ast \Gr^i \cV
$$
of $p^\ast \cV$. So it suffices to trivialize each $p^\ast \Gr^i \cV$.

To do this, we first do it on a single leaf $\L$ of $P$. The
restriction of the monodromy representation $\tau$ to $\L$
is clearly trivial. Consequently, the restriction of $p^\ast \cV$
to $\L$ is trivial as a flat bundle. Observe that if this leaf contains
$\xtilde_o$, then this trivialization satisfies  condition (\ref{cond})
in the statement of (\ref{triv}).

Next, change the trivialization of $p^\ast\Gr^i\cV$ on $p^{-1}(x_o)$
so that it satisfies condition (\ref{cond}) in the statement of
(\ref{triv}). Extend this to a trivialization of $p^{-1}\Gr^i\cV$
on all of $P$ by parallel transport along the leaves of $P$. This
gives a well defined local trivialization which is a global
trivialization by the argument in the previous paragraph.

We thus obtain a trivialization of $p^\ast\cV$ which is
compatible with the filtration $\cV^\dt$ and which is flat on
each $\Gr^i\cV$. It follows that the connection form $\omegatilde$
associated to this trivialization satisfies
$\omegatilde \in E^1(P)\otimes \u$.
\end{proof}

If $S$ is not finite, this connection is not flat as it is
not flat in the vertical direction. We can make it flat by
restricting it to the leaves of the foliation $\F$ of $P$.
Denote the image of $\omegatilde$ under the restriction
homomorphism
$$
E^1(P)\otimes \u \to E^1(P,\F) \otimes \u
$$
by $\omega$. It defines the connection in the leaf direction.
This connection is clearly flat, and it follows that $\omega$ is
integrable.

The following assertion is a consequence of the
properties (1) and (2) in the statement of Proposition \ref{triv}
and (\ref{dga_homom}). Note that we view $S$ as acting on the
left of $\u$ via the adjoint action --- that is, via the
composite $S \to \im \tau \hookrightarrow G \to \Aut \u$.

\begin{proposition}
The connection form $\omega$ is integrable and lies in the subspace
$\Efin^1(X,\O(P))\otimes \u$ of $E^1(P,\F)\otimes \u$. Moreover,
if $s\in S$, then
$s^\ast \omega = Ad(s)\omega$. \qed
\end{proposition}

\begin{remark}\label{converse}
There is a converse to this result. Suppose that $\u$ is a nilpotent
Lie algebra in the category of rational representations of $S$. Then
we can form the semi-direct product $G=S\ltimes U$, where $U$ is the
corresponding unipotent group. If $V$ is a $G$ module, and if
$$
\omega \in \Efin^1(X,\O(P))\otimes \u
$$
satisfies the conditions
\begin{enumerate}
\item $d\omega + \omega \wedge \omega = 0$;
\item $s^\ast \omega = Ad(s)\omega$;
\end{enumerate}
then we can construct an object of $\B(X,S)$ with fiber $V$ over $x_o$
whose pullback to $P$ has connection form $\omega$ with respect to an
appropriate trivialization.
\end{remark}

We are now ready to express the monodromy representation of
$\V$ in terms of iterated integrals of $\omega$. Recall that
K.-T.~Chen \cite{chen} defined, for 1-forms $w_i$ on a manifold $M$
taking values in an associative algebra $A$,
$$
\int_\gamma w_1 w_2 \dots w_r
$$
to be the element
$$
\idotsint\limits_{0 \le t_1 \le \dots \le t_r \le 1}
f_1(t_1)f_2(t_2) \dots f_r(t_r)\, dt_1dt_2 \dots dt_r
$$
of $A$. This is regarded as an $A$-valued function $PM \to A$ on the
path space of $M$. An $A$-valued iterated integral is a function
$PM \to A$ which is a linear combination of functions
of this form together with a constant function.

Suppose that $V\times M \to M$ is a trivial bundle with a connection
given by the connection form
$$
\omega \in E^1(M) \otimes \End(V).
$$
In this case we can define the parallel transport map
$$
T : PM \to \Aut(V)
$$
where $PM$ denotes the space of piecewise smooth paths in $M$. A
path goes to the linear transformation of $V$ obtained by parallel
transporting the identity along it. Chen \cite{chen} obtained the
following expression for $T$ in terms of $\omega$.

\begin{proposition}\label{transp}
With notation as above, we have
$$
T(\gamma) = 1 +  \int_\gamma \omega + \int_\gamma \omega\omega
+ \int_\gamma \omega\omega\omega + \cdots \qed
$$
\end{proposition}

Note that since $\u$ is nilpotent, this is a finite sum.
Armed with this formula, we can express the monodromy
of $\V \to X$ in terms of $\omega \in \Efin^1(X,\O(P))$. Suppose
that $\gamma \in P_{x_o,x_o}X$. Denote the unique lift of
$\gamma$ to $P$ which is tangent to  $\F$ and begins at
$\xtilde_o\in p^{-1}(x_o)$, by $\gammatilde$.

\begin{proposition}\label{monod}
The monodromy of $\V \to X$ takes $\gamma \in P_{x_o,x_o}X$
to
$$
\rhotilde(\gamma) =
\left(1 + \int_{\gammatilde} \omega + \int_{\gammatilde}
\omega\omega + \int_{\gammatilde} \omega\omega\omega + \cdots
\right)\tau(\rho(\gamma)) \in G.
$$
\end{proposition}

The proof is a straightforward consequence of Chen's formula
(\ref{transp}) and condition (\ref{cond}) of (\ref{triv}).

This formula motivates the following generalization of Chen's iterated
integrals.

\begin{definition}\label{defn}
For $\phi\in \O(S)$ and $w_1,\dots,w_r$ elements of
$\Efin^1(X,\O(P))$, we define
$$
\int \left(w_1 \dots w_r | \phi\right) : P_{x_o,x_o}X \to \R
$$
by
$$
\int_\gamma \left(w_1 \dots w_r | \phi\right)
= \phi(\rho(\gamma))\int_{\gammatilde}w_1\dots w_r.
$$
\end{definition}

We will call linear combinations of such functions {\it iterated
integrals with coefficients in $\O(S)$}. They will be regarded
as functions $P_{x_o,x_o}X \to \R$. We will denote the set of them
by $I(X,\O(S))_{x_o}$. Such an iterated integral will
be said to be {\it locally constant} if it is constant on each
connected component of $P_{x_o,x_o}X$. We shall denote the set of
locally constant iterated integrals on $P_{x_o,x_o}X$ by
$H^0(I(X,\O(S))_{x_o})$. Evidently, each such locally
constant iterated iterated integral defines a function
$$
\pi_1(X,x_o) \to \R.
$$
By taking matrix entries in (\ref{monod}), we obtain the following
result.

\begin{corollary}\label{matrix_entries}
Each matrix entry of the monodromy representation
$$
\rhotilde: \pi_1(X,x_o) \to G
$$
of an object of $\B(X,S)$ can be expressed as a locally constant
iterated integral on $X$ with coefficients in $\O(S)$. \qed
\end{corollary}

The following results imply that $H^0(I(X,\O(S))_{x_o})$
is a Hopf algebra with coproduct dual to the multiplication of
paths, and antipode dual to the involution of $P_{x_o,x_o}X$ that
takes each path to its inverse.

\begin{proposition}\label{props}
Suppose that $\gamma$ and $\mu$ are in $P_{x_o,x_o}X$, that
$\phi,\psi\in \O(S)$ and that $w_1,w_2,\dots \in \Efin^1(X,\O(P))$.
Then we have:
\begin{equation}
\int_\gamma \left(w_1 \dots w_p | \phi\right)
\int_\gamma \left(w_{p+1} \dots w_{p+q} | \phi\right)
= \sum_{\sigma \in Sh(p,q)}
\int_\gamma \left( w_{\sigma(1)}\dots w_{\sigma(p+q)}|
\phi\psi \right)
\end{equation}
where $Sh(p,q)$ denotes the set of shuffles of type $(p,q)$;
\begin{equation}
\int_{\gamma^{-1}} \left(w_1 \dots w_r | \phi\right) = (-1)^r
\int_\gamma
\left(\rho(\gamma^{-1})^\ast w_r \dots \rho(\gamma^{-1})^\ast w_1|
i_S^\ast\phi\right)
\end{equation}
where $i_S^\ast : \O(S) \to \O(S)$ is the antipode of $\O(S)$;
\begin{equation}
\int_{\gamma\mu}\left(w_1 \dots w_r | \phi\right) =
\sum_{i=0}^r \sum_j \int_\gamma\left( w_1\dots w_i| \phi_j'\right)
\int_\mu \left(\rho(\gamma)^\ast w_{i+1}\dots \rho(\gamma)^\ast w_r |
\phi_j''\right)
\end{equation}
where $\Delta_S : \O(S) \to \O(S)\otimes \O(S)$ is the coproduct of
$\O(S)$, and
$$
\Delta_S \phi = \sum_j \phi_j'\otimes \phi_j''.
$$
\end{proposition}

\begin{proof}
This proof is a straightforward using the definition (\ref{defn})
and basic properties of classical iterated integrals due to Chen
\cite{chen}.
\end{proof}

\begin{corollary}
The set of iterated integrals $I(X,\O(S))_{x_o}$ is a commutative Hopf
algebra.
\end{corollary}

\begin{remark}
Let $\pi_1(X,x_o) \to \G$ be the completion of $\pi_1(X,x_o)$ relative
to $\rho : \pi_1(X,x_o) \to S$. Since the coordinate ring $\O(\G)$ of
$\G$ is the ring of matrix entries of representations of $G$, it follows
from (\ref{matrix_entries}) that there is a Hopf algebra inclusion
$$
\O(\G) \hookrightarrow H^0(I(X,\O(S))_{x_o}).
$$
To prove this assertion, it would suffice to show that
$$
H^0(I(X,\O(P)))\otimes_{\O(S)}\R
$$
is the direct limit of coordinate rings of a directed system of unipotent
groups, each with an $S$ action. This is surely true, but we seek a more
algebraic de~Rham theorem for $\O(\G)$ which is more convenient for Hodge
theory.
\end{remark}

\section{Higher Iterated Integrals}

As a preliminary step to defining the algebraic analogue of
$I(X,\O(S))_{x_o}$, we generalize the the definition of iterated
integrals with values in $\O(S)$ to higher dimensional forms.

Denote by $E^n(P_{x_o,x_o}X)$ the differential forms of degree $n$ on
the loop space $P_{x_o,x_o}X$. One can surely use any reasonable
definition of differential forms on $P_{x_o,x_o}X$, but we will use
Chen's definition from \cite{chen} where, to specify a differential form
on $P_{x_o,x_o}X$, it is enough to specify its pullback along each
``smooth map'' $\alpha : U\to P_{x_o,x_o}X$ from an open subset $U$
of some finite dimensional euclidean space. By a smooth map,
we mean a map $\alpha : U\to P_{x_o,x_o}X$ whose ``suspension''
$$
\alphahat : [0,1] \times U \to X; \quad (t,u) \mapsto \alpha(u)(t)
$$
is continuous and smooth on each $[t_{j-1},t_j]\times U$ for some
partition
$$
0=t_0 \le t_1 \le \dots \le t_m = 1
$$
of $[0,1]$.

\begin{definition}\label{higherdef}
Suppose that $\phi \in \O(S)$, and that $w_j\in \Efin^{n_j}(X,\O(P))$ with
each $n_j>0$. Set $n=-r + \sum_j n_j$. Define
$$
\int \left(w_1 \dots w_r | \phi\right) \in E^n(P_{x_o,x_o}X)
$$
by specifying that for each smooth map $\alpha : U \to P_{x,x}X$,
$$
\alpha^\ast\int \left(w_1 \dots w_r | \phi\right)
$$
is the element
$$
\idotsint\limits_{0\le t_1 \le \cdots \le t_r \le 1}
\what_1(t_1)\wedge \dots \wedge \what_r(t_r)\,
dt_1\! dt_2\dots dt_r\, \phi(\rho(\alpha(u)))
$$
of $E^n(U)$, where
$$
\what_j : (\partial/\partial t) \lrcorner \alphatilde^\ast w_j
$$
and $\alphatilde : [0,1]\times U \to P$ is the smooth map with the
property that for each $x\in U$, the map $t\mapsto \alphatilde(t,x)$
is the unique lift of $t\mapsto \alphahat(t,x)$ that begins at
$\xtilde_o$ and is tangent to $\F$.
\end{definition}

These iterated integrals form a subspace $I^\dt(X,\O(S))_{x_o}$
of $E^\dt(P_{x_o,x_o}X)$. Chen's arguments \cite{chen} can be adapted
easily to show that this is, in fact, a sub d.g.\ Hopf algebra of
$E^\dt(P_{x_o,x_o}X)$. In particular, we have:

\begin{proposition}
The space of locally constant iterated integrals on $X$ with
coefficients in $\O(S)$ is $H^0(I^\dt(X,\O(S))_{x_o})$. \qed
\end{proposition}

\section{The Reduced Bar Construction}

In this section we review Chen's definition of the reduced bar
construction which he described in \cite{chen:bar}.

Suppose that $A^\dt$ is a commutative differential graded algebra
(hereafter denoted d.g.a.) and that $M^\dt$ and $N^\dt$ are
complexes which are modules over $A^\dt$. That is, the
structure maps
$$
A^\dt \otimes M^\dt \to M^\dt \text{ and }
A^\dt\otimes N^\dt \to N^\dt
$$
are chain maps. We shall suppose that $A^\dt$, $M^\dt$ and $N^\dt$
are all positively graded. Denote the subcomplex of $A^\dt$ consisting
of elements of positive degree by $A^+$.

The {\it (reduced) bar construction} $B(M^\dt,A^\dt,N^\dt)$ is
defined as follows. We first describe the underlying graded vector
space. It is a quotient of the graded vector space
$$
T(M^\dt,A^\dt,N^\dt) :=
\bigoplus_s M^\dt \otimes\left(A^+[1]^{\otimes r}\right) \otimes N^\dt.
$$
We will use the customary notation $m[a_1|\dots|a_r]n$ for
$$
m\otimes a_1\otimes \dots \otimes a_r \otimes n
\in T(M^\dt,A^\dt,N^\dt).
$$
To obtain the vector space underlying the bar construction, we mod out
by the relations
$$
m[dg|a_1|\dots|a_r]n = m[ga_1|\dots|a_r]n - m\cdot g[a_1|\dots|a_r]n;
$$
\begin{multline*}
m[a_1|\dots|a_i|dg|a_{i+1}|\dots|a_r]n = \hfill \cr
m[a_1|\dots|a_i|g\,a_{i+1}|\dots|a_r]n
- m[a_1|\dots|a_i\,g|a_{i+1}|\dots|a_r]n \quad 1\le i < s;
\end{multline*}
$$
m[a_1|\dots|a_r|dg]n = m[a_1|\dots|a_r]g\cdot n - m[a_1|\dots|a_r\,g]n;
$$
$$
m[dg]n = 1 \otimes g\cdot n  - m\cdot g \otimes 1
$$
Here each $a_i \in A^+$, $g\in A^0$, $m\in M^\dt$,  $n\in N^\dt$,
and $r$ is a positive integer.

Before defining the differential, it is convenient to define an
endomorphism $J$ of each graded vector space by
$J: v\mapsto (-1)^{\deg v}v$. The differential is defined as
$$
d = d_M\otimes 1_T \otimes 1_N + J\otimes d_B \otimes 1 +
J_M \otimes J_T \otimes d_N + d_C.
$$
Here $T$ denotes the tensor algebra on $A^+[1]$, $d_B$ is defined by
\begin{multline*}
d_B[a_1|\dots|a_r] =
\sum_{1\le i \le r} (-1)^i [Ja_1|\dots|Ja_{i-1}|da_i|a_{i+1}|\dots|a_r]
\cr \hfill + \sum_{1 \le i < r}
(-1)^{i+1}[Ja_1|\dots|Ja_{i-1}|Ja_i\wedge a_{i+1}|a_{i+2}|\dots|a_r]
\end{multline*}
and $d_C$ is defined by
$$
d_C m[a_1|\dots|a_r]n = (-1)^s
Jm[Ja_1|\dots|Ja_{r-1}]a_r \cdot n - Jm\cdot a_1 [a_2|\dots|a_r]n.
$$
One can check that these differentials are well defined.

If both $M^\dt$ and $N^\dt$ are themselves d.g.a.s over $A^\dt$,
then $B(M^\dt,A^\dt,N^\dt)$ is also a differential graded algebra.
The product is defined by
\begin{multline}\label{prod}
m'[a_1|\dots|a_p]n' \otimes m''[a_{p+1}|\dots |a_{p+q}]n'' \mapsto
\hfill \cr \hfill \sum_{\sigma\in \Sigma(p,q)} \pm m'\wedge m''
[a_{\sigma(1)}| a_{\sigma(2)}|\dots |a_{\sigma(p+q)}]n'\wedge n''.
\end{multline}
Here $\Sigma(p,q)$ denotes the set of shuffles of type $(p,q)$. The
sign in front of each term on the right hand side is determined by
the usual sign conventions that apply when moving a symbol of degree
$k$ past one of degree $l$ --- one considers each $a_j$ to be of
degree $-1 + \deg a_j$.

The reduced bar construction $B(M^\dt,A^\dt,N^\dt)$ has a standard
filtration
$$
\R = B_0(M^\dt,A^\dt,N^\dt) \subseteq B_1(M^\dt,A^\dt,N^\dt)
\subseteq B_2(M^\dt,A^\dt,N^\dt) \subseteq \cdots
$$
which is often called the {\it bar filtration}. The subspace
$$
B_s(M^\dt,A^\dt,N^\dt)
$$
is defined to be the span of those $m[a_1|\dots|a_r]n$ with $r\le s$.
When $A^\dt$ has connected homology (i.e., $H^0(A^\dt) = \R$), the
corresponding spectral sequence, which is called the {\it Eilenberg-Moore
spectral sequence}, has $E_1$ term
$$
E_1^{-s,t} =
\left[M^\dt\otimes H^+(A^\dt)^{\otimes s}\otimes N^\dt\right]^t.
$$
A proof of this can be found in \cite{chen:bar}.

The following basic property of the reduced bar construction is
a special case of a result proved in \cite{chen:bar}. It is easily
proved using the Eilenberg-Moore spectral sequence. Suppose that
$\psi : A_1^\dt \to A_2^\dt$ is a d.g.a.\ homomorphism, and that
$M^\dt$ is a right $A_2^\dt$ module and $N^\dt$ a right $A_2^\dt$
module. Then $M^\dt$ and $N^\dt$ can be regarded as $A_1^\dt$ modules
via $\psi$. We therefore have a chain map
\begin{equation}\label{map}
B(M^\dt,A_1^\dt,N^\dt) \to B(M^\dt,A_2^\dt,N^\dt).
\end{equation}

\begin{proposition}\label{qism}
If $\psi$ is a quasi-isomorphism, then so is (\ref{map}). \qed
\end{proposition}

\section{The Construction of $\Gdr$}
\label{gdr}

In this section we construct a proalgebraic group $\Gdr$ which
is an extension
$$
1 \to \Udr \to \Gdr \stackrel{p}{\to} S \to 1,
$$
where $\Udr$ is prounipotent, and a homomorphism
$$
\rhotilde : \pi_1(X,x_o) \to \Gdr
$$
whose composition with $p : \Gdr \to S$ is $\rho$. We do this
by constructing the coordinate ring of $\Gdr$ using the
bar construction. In the two subsequent sections, we will show
that $\rhotilde : \pi_1(X,x_o) \to \Gdr$
is the Malcev completion of $\pi_1(X,x_o)$ relative to $\rho$.

The fixed choice of a base point $\xtilde_o \in p^{-1}(x_o)$
determines augmentations
$$
\epsilon_{\xtilde_o} : \Efin^\dt(X,\O(P)) \to \R
$$
and
$$
\delta_{x_o} : \Efin^\dt(X,\O(P)) \to \O(S)
$$
as we shall now explain.
Since these augmentations are compatible with restriction, it suffices
to give them in a neighbourhood of $x_o$. Over a contractible
neighbourhood $U$ of $x_o$, the local system $P$ is trivial and
may therefore by identified with the trivial flat bundle
$S\times U \to S$ in such a way that $\xtilde_o$ corresponds to
$(1,x_o) \in S\times U$. The restriction of an element of
$\Efin^k(X,\O(P))$ to $U$ is then of the form
$$
\sum_i \phi_i \otimes w_i
$$
where $\phi_i \in \O(S)$, and $w_i \in E^k(U)$.  Denote the augmentation
$E^\dt(U) \to \R$ induced by $x_o$ by $\mu_{x_o}$. Then the augmentations
$\delta_{x_o}$ and $\epsilon_{\xtilde_o}$ are defined by
$$
\delta_{x_o} : \sum_i \phi_i \otimes w_i
\mapsto \sum_i \mu_{x_o}(w_i)\, \phi_i
$$
and
$$
\epsilon_{\xtilde_o} : \sum_i \phi_i \otimes w_i \mapsto
\sum_i \mu_{x_o}(w_i)\, \phi_i(1).
$$

One can regard $\R$ and $\O(S)$ as algebras over $\Efin^\dt(X,\O(P))$
where the actions of $\Efin^\dt(X,\O(P))$ on these is defined using
these two augmentations. We can therefore form
the bar construction
$$
B(\R,\Efin^\dt(X,\O(P)),\O(S))
$$
which we shall denote by $B(\Efin^\dt(X,\O(P))_{\xtilde_o,(x_o)})$.
It is a commutative d.g.a.\ when endowed with the product (\ref{prod}).
It is, in fact, a d.g.~Hopf algebra, with coproduct defined as follows:
\begin{multline*}
\Delta : [w_1|\dots|w_r]\phi \mapsto\\
\sum_i [w_1|\dots |w_i]
\left(\sum \psi_i^{(k_i)}\dots\psi_r^{(k_r)}\phi'\right)
\otimes [w_{i+1}^{(k_i)}|\dots|w_r^{(k_r)}]\left(\sum \phi''\right)
\end{multline*}
where
$$
\Delta_S(\phi) = \sum \phi'\otimes \phi''
$$
is the diagonal of $\O(S)$, and the map
$$
\Efin^\dt(X,\O(P)) \to \O(S) \otimes \Efin^1(X,\O(P))
$$
which gives the $S$ action takes $w_j$ to
$$
\sum \psi_j^{(k_j)} \otimes w_j^{(k_j)}.
$$

The following proposition is a direct consequence of the definition
(\ref{higherdef}) and the basic properties of iterated integrals
which may be found in \cite{chen}.

\begin{proposition}\label{pro}
The map
$$
B(\Efin^\dt(X,\O(P))_{\xtilde_o,(x_o)}) \to I^\dt(X,\O(S))_{x_o}
$$
defined by
$$
[w_1|w_2|\dots |w_r]\phi \mapsto \int\left(w_1 w_2 \dots w_r|\phi\right)
$$
is a well defined d.g.~Hopf algebra homomorphism. \qed
\end{proposition}

\begin{proposition}\label{alg_gp}
If $\pi_1(X,x_o)$ is finitely generated, then
$$
H^0(B(\Efin^\dt(X,\O(P))_{\xtilde_o,(x_o)}))
$$
is the coordinate ring of a linear proalgebraic group which is an
extension of $S$ by a prounipotent group.
\end{proposition}

In the proof, we shall need the following technical result, the proof of
which is a straightforward modification of Sullivan's proof of the
existence of minimal models (cf.\ \cite{sullivan}.)

\begin{proposition}\label{sub}
There is a d.g.~subalgebra $A^\dt$ of $\Efin^\dt(X,\O(P))$ with
$A^0=\R$, which is also an $S$ submodule, with the properties that the
inclusion is a quasi-isomorphism. \qed
\end{proposition}

\begin{proof}[Proof of (\ref{alg_gp})]
Choose a d.g.\ subalgebra $A^\dt$ of $\Efin^\dt(X,\O(P))$ as given by
(\ref{sub}). It follows from (\ref{qism}) that the natural map
$$
H^0(B(\R,A^\dt,\O(S))) \to H^0(B(\Efin^\dt(X,\O(P))_{\xtilde_o,(x_o)}))
$$
is an isomorphism. Since $A^0=\R$, we have that
$$
H^0(B(\R,A^\dt,\O(S))) = H^0(B(\R,A^\dt,\R))\otimes \O(S).
$$
It is not difficult to check that $\O(S)$ is a sub Hopf algebra, and
that this is a tensor product of algebras, but where the coproduct is
twisted by the action of $S$ on $H^0(B(\R,A^\dt,\R))$. So, if we can
show that $H^0(B(\R,A^\dt,\R))$ is the limit of the coordinate rings
of an inverse system of unipotent groups, each with an $S$ action, then
we will have shown that
$$
H^0(B(\Efin^\dt(X,\O(P))_{\xtilde_o,(x_o)}))
$$
is the coordinate ring of
$$
S\ltimes \spec H^0(B(\R,A^\dt,\R))
$$
and therefore proved the proposition. From \cite{hain:bar}, we
know that there is a canonical splitting (in particular, it is
$S$ equivariant) of the projection
$$
H^0(B(\R,A^\dt,\R)) \to QH^0(B(\R,A^\dt,\R)) =: Q
$$
onto the indecomposable elements $Q$. This splitting induces an
$S$-equivariant algebra isomorphism
$$
\R[Q] \to H^0(B(\R,A^\dt,\R)).
$$
The bar filtration induces a filtration
$$
Q_1\subseteq Q_2 \subseteq Q_3 \subseteq \cdots \subseteq Q
$$
of the indecomposables such that $Q= \cup Q_r$. Each $Q_r$ is a Lie
coalgebra, and the cobracket $\Delta^c$ satisfies
$$
\Delta^c : Q_r \to \sum_{i+j = r} Q_i \otimes Q_j
$$
and is injective when $r>1$. Since $\pi_1(X,x_o)$ is finitely generated,
each of the cohomology groups
$H^1(X,\V)$ is finite dimensional for each rational local system $\V$
over $X$. It follows from (\ref{isom}) that each isotypical component of
$H^1(X,\O(P))$ is finite dimensional. Since $Q_r/Q_{r-1}$ is a
subquotient of
$$
H^1(X,\O(P))^{\otimes r},
$$
it follows that each $S$-isotypical component of each $Q_r$ is finite
dimensional. One can now prove by induction on $r$ using the nilpotence,
that as an $S$-module, each $Q_r$ is the direct limit of duals of
nilpotent Lie algebras, each of which has an $S$ action. This completes
the proof.
\end{proof}

\begin{definition}\label{def}
Define proalgebraic groups $\Gdr$ and $\Udr$ by
$$
\Gdr = \spec H^0(B(\Efin^\dt(X,\O(P))_{\xtilde_o,(x_o)}))
$$
and
$$
\Udr = \spec H^0(B(\R,\Efin^\dt(X,\O(P)),\R)).
$$
\end{definition}

Evidently, we have an extension
$$
1 \to \Udr \to \Gdr \to S \to 1
$$
of proalgebraic groups, where $\Udr$ is prounipotent.

When we want to emphasize the dependence of $\Gdr$ and $\Udr$
on $(X,x)$, we will write them as $\Gdr(X,x)$ and $\Udr(X,x)$,
respectively.

\begin{proposition}\label{homom}
There is a natural homomorphism $\rhotilde : \pi_1(X,x_o) \to \Gdr$ whose
composition with $\Gdr \to S$ is $\rho$.
\end{proposition}

\begin{proof}
Define a map from $P_{x_o,x_o}X$ to the linear functionals on
$$
B(\Efin^\dt(X,\O(P))_{\xtilde_o,(x_o)})
$$
by
$$
\gamma : [w_1|\dots|w_r]\phi
\mapsto \int_\gamma\left(w_1\dots w_r|\phi\right).
$$
This induces a function
$$
\Phi : \pi_1(X,x_o) \to \Hom_\R(\O(\Gdr),\R).
$$
Define $\rhotilde$ by taking the class of $\gamma$ in $\pi_1(X,x_o)$
to the maximal ideal of
$$
H^0(B(\Efin^\dt(X,\O(P))_{\xtilde_o,(x_o)}))
$$
consisting of those elements on which $\gamma$ vanishes. (Note
that $\gamma$ acts via integration.) That this
is a group homomorphism follows from (\ref{props}).
\end{proof}

\section{Construction of Homomorphisms from $\Gdr$}

Suppose that $G$ is a linear algebraic group which can be
expressed as an extension
$$
1 \to U \to G \to S \to 1
$$
where $U$ is unipotent. Choose an isomorphism of $G$ with $S\ltimes U$.
Denote the Lie algebra of $U$ by $\u$.

\begin{proposition}
Each one form $\omega \in \Efin^1(X,\O(P))\otimes \u$ that satisfies
\begin{enumerate}
\item $d\omega + \omega\wedge \omega = 0$;
\item for all $s \in S$, $s^\ast \omega = Ad(s)\omega$;
\end{enumerate}
determines a homomorphism $\Gdr \to G$ that commutes with projection to
$S$.
\end{proposition}

\begin{proof}
First note that since the exponential map $\u \to U$ is a polynomial
isomorphism, $\O(U)$ is isomorphic to the polynomials $\R[\u]$ on the
vector space $\u$. Further, there is a natural isomorphism
\begin{equation}\label{iso}
\O(U) \cong \R[\u] \to \lim_\to \Hom(U\u/I^n,\R)
\end{equation}
which is defined by noting that $U\u$ is, by the PBW Theorem, the
symmetric coalgebra $S^c\u$ on $\u$. The isomorphism (\ref{iso})
is an isomorphism of Hopf algebras.

Set
$$
T = 1 + [\omega] + [\omega|\omega] + [\omega|\omega|\omega] + \cdots
$$
which we view as an element of
$$
B(\Efin^\dt(X,\O(P))_{\xtilde_o,(x_o)}) \comptensor \Uhat\u
$$
of degree zero, where $\Uhat\u$ denotes the completion
$$
\lim_\leftarrow U\u/I^n
$$
of $U\u$ with respect to the powers of its augmentation ideal,
and where $\comptensor$ denotes the completed tensor product
$$
\lim_\leftarrow
B(\Efin^\dt(X,\O(P))_{\xtilde_o,(x_o)})\otimes \Uhat\u/I^n
$$

The coordinate ring of $G$ is isomorphic to $\O(U)\otimes\O(S)$.
Define a linear map
$$
\Theta : \O(G) \to B(\Efin^\dt(X,\O(P))_{\xtilde_o,(x_o)})^0
$$
by
$$
f\otimes \phi \mapsto \langle T,f\rangle \cdot \phi.
$$
It is not difficult to check that $\Theta$ is a well defined Hopf algebra
homomorphism. This uses the fact that $s^\ast \omega = Ad(s)\omega$. That
$\omega$ satisfies the integrability condition
$$
d \omega + \omega \wedge \omega = 0
$$
implies that $\im \Theta$ is contained in
$H^0(B(\Efin^\dt(X,\O(P))_{\xtilde_o,(x_o)}))$. It follows that $\Theta$
induces a Hopf algebra homomorphism
$$
\O(G) \to H^0(B(\Efin^\dt(X,\O(P))_{\xtilde_o,(x_o)}))
$$
and therefore a group homomorphism
$$
\theta : \Gdr \to G
$$
which commutes with the projections to $S$.

Finally, it follows from (\ref{converse}), (\ref{transp}) and (\ref{homom})
that the composite
$$
\pi_1(X,x_o) \to \Gdr \to G
$$
is the homomorphism induced by $\omega$.
\end{proof}

\begin{corollary}\label{monod_rep}
If $\V$ is a local system in $\B(X,S)$, then the monodromy representation
$$
\tau : \pi_1(X,x_o) \to \Aut V_o
$$
factors through $\rhotilde : \pi_1(X,x_o) \to \Gdr$. \qed
\end{corollary}

\begin{corollary}\label{factor}
If $\tau: \pi_1(X,x_o) \to G$ is a homomorphism into a linear algebraic
group which is an extension of $S$ by a unipotent group, and whose
composite with the projection to $S$ is $\rho$, then there is a
homomorphism $\Gdr \to G$ whose composite with
$$
\rhotilde : \pi_1(X,x_o) \to \Gdr
$$
is $\tau$.
\end{corollary}

\begin{proof}
Denote the kernel of $G \to S$ by $U$. One can construct a faithful, finite
dimensional representation $V$ of $G$ which has a filtration
$$
V = V^0 \supset V^1 \supset V^2 \supset \cdots
$$
by $G$-submodules whose intersection is zero and where each $V^j/V^{j+1}$
is a trivial $U$-module. The corresponding local system over $X$ lies in
$\B(X,S)$. The result now follows from (\ref{monod_rep}).
\end{proof}

\section{Isomorphism with the Relative Completion}

Denote $\pi_1(X,x_o)$ by $\pi$. In Section~\ref{gdr} we constructed
a homomorphism $\pi \to \Gdr$. In this section, we prove:

\begin{theorem}\label{main}
If $\pi$ is finitely generated, then the homomorphism $\pi \to \Gdr$
is the completion of $\pi$ relative to $\rho$.
\end{theorem}

To prove the theorem, we first fix a completion
$\pi \to \G$ of $\pi$ relative to $\rho$.
The universal mapping property of the relative completion gives a
homomorphism $\G \to \Gdr$ of proalgebraic groups that commutes with
the canonical projections to $S$. It follows from  (\ref{factor}), that
there is a natural homomorphism
\begin{equation}\label{univ}
\Gdr \to \G
\end{equation}
that also commutes with the projections to $S$. It follows from the
universal mapping property of the relative completion that the composite
$$
\G \to \Gdr \to \G
$$
is the identity.

Denote the prounipotent radical of $\G$ by $\U$. Since $\pi$ is finitely
generated, each of the groups $H^1(\pi,V)$ is finite dimensional for
each rational representation $V$ of $S$.

In view of the following proposition and the assumption that $\pi$
is finitely generated, all we need do to show that the natural
homomorphism $\G \to \Gdr$ is an isomorphism is to show that either
of the induced maps
$$
\Hom_S(H_1(\U),V) \to \Hom_S(H_1(\Udr),V) \to \Hom_S(H_1(\U),V)
$$
is an isomorphism for all $S$-modules $V$.

\begin{proposition}
Suppose that $G_1$ and $G_2$ are extensions of the reductive group $S$ by
unipotent groups $U_1$, $U_2$, respectively:
$$
1 \to U_j \to G_j \to S \to 1.
$$
Suppose that $\theta : G_1 \to G_2$ is a split surjective homomorphism of
algebraic groups that commutes with the projections to $S$. If either of
the induced maps
$$
\Hom_S(H_1(U_1),V) \to \Hom_S(H_1(U_2),V) \to \Hom_S(H_1(U_1),V)
$$
is an isomorphism for all $S$ modules $V$, then both are, and $\theta$
is an isomorphism.
\end{proposition}

\begin{proof}
The proof reduces to basic fact that a split surjective homomorphism
between nilpotent Lie algebras is an isomorphism if and only if it
induces an isomorphism on $H_1$. The details are left to the reader.
\end{proof}

Our first task in the proof of Theorem~\ref{main} is to
compute $\Hom_S(H_1(\U),V)$.

\begin{proposition}\label{h1_comp}
For all $S$-modules $V$, there is a canonical isomorphism
$$
H^1(\pi,V) \cong \Hom_S(H_1(\U),V).
$$
\end{proposition}

\begin{proof}
We introduce an auxiliary group for the proof. Let
$$
\Hom_\rho(\pi,S\ltimes V)
$$
be the set of group homomorphisms $\pi \to S\ltimes V$ whose composite
with the projection $S\ltimes V \to S$ is $\rho$. Then there is a natural
bijection between  $\Hom_\rho(\pi,S\ltimes V)$ and the set of splittings
$\pi \to \pi \ltimes V$ of the projection $\pi \ltimes V \to \pi$: the
splitting $\sigma$ corresponds to $\rhotilde : \pi \to S\ltimes V$ if and
only if the diagram
$$
\begin{CD}
\pi  @>\sigma>> \pi\ltimes V \cr
@| @VV{\rho\ltimes id}V\cr
\pi @>\rhotilde>> S\ltimes V \cr
\end{CD}
$$
commutes.

The kernel $V$ acts on both $\Hom_\rho(\pi,S\ltimes V)$ and on the set of
splittings, in both cases by inner automorphisms. The action commutes
with the bijection. Since $H^1(\pi,V)$ is naturally isomorphic to the set
of splittings of $\pi\ltimes V \to \pi$ modulo conjugation by $V$
\cite[p.~106]{maclane}, the bijection induces a natural isomorphism
$$
H^1(\pi,V) \cong \Hom_\rho(\pi,S\ltimes V)/\sim.
$$

On the other hand, by the universal mapping property of the relative
completion, each element of $\Hom_\rho(\pi,S\ltimes V)$ induces a
homomorphism $\G \to S\ltimes V$ which commutes with the projections to
$S$. Such a homomorphism induces a homomorphism $\U \to V$, and therefore an
$S$-equivariant homomorphism $H_1(\U) \to V$. Since $V$ is central, this
induces a homomorphism
$$
\Hom_\rho(\pi,S\ltimes V) \to \Hom_S(H_1(\U),V).
$$
To complete the proof, we show that this is an isomorphism. Denote the
commutator subgroup of $\U$ by $\U'$. Then the quotient $\G/\U'$ is an
extension of $S$ by $H_1(\U)$; the latter being a possibly infinite
product of representations of $S$ in which each isotypical factor is
finite dimensional. Using the fact that every extension of $S$ by a
rational representation in the category of algebraic groups splits and
that any two such splittings are conjugate by an element of the kernel,
we see that the extension
\begin{equation}\label{seqce}
0 \to H_1(\U) \to \G/\U' \to S \to 1
\end{equation}
is split and that any two splittings are conjugate by an element of
$H_1(\U)$. Choose a splitting of this sequence. This gives an isomorphism
$$
\G/\U' \cong S\ltimes H_1(\U).
$$
An $S$-equivariant homomorphism $H_1(\U) \to V$ induces a homomorphism
$$
\G/\U' \cong S\ltimes H_1(\U) \to S\ltimes V
$$
of proalgebraic groups. Composing this with the homomorphism
$$
\pi \to \G \to \G/\U',
$$
we obtain an element of $\Hom_\rho(\pi,S\ltimes V)/\sim$. Since all
splittings of (\ref{seqce}) differ by an inner automorphism by an element
of $H_1(\U)$, we have constructed a well defined map
$$
\Hom_S(H_1(\U),V) \to \Hom_\rho(\pi,S\ltimes V).
$$
This is easily seen to be the inverse of the map constructed above. This
completes the proof.
\end{proof}

The following result completes the proof of Theorem~\ref{main}.

\begin{proposition}
The map
$$
H_1(\U) \to H_1(\Udr)
$$
induced by (\ref{univ}) is an isomorphism.
\end{proposition}

\begin{proof}
It suffices to show that for all rational representations $V$ of
$S$, the map
$$
\left[H^1(\Udr)\otimes V\right]^S \to \left[H^1(\U)\otimes V\right]^S
$$
is an isomorphism. Both groups are isomorphic to $H^1(X,\V)$.

Choose a de~Rham representative $w \in E^1(X,\V)$ of a class in
$H^1(X,\V)$. Let $\delta \in V^\ast\otimes V$ be the element
corresponding to the identity $V\to V$. Set
$$
\omega := w\otimes \delta \in E^1(X,\V)\otimes V^\ast \otimes V.
$$
Regard $V$ as an abelian Lie algebra. Then
$$
\omega \in \Efin^1(X,\O(P))\otimes V.
$$
It is closed, and therefore satisfies the integrability condition
$d\omega + \omega\wedge\omega = 0$. Since the identity $V\to V$ is
$S$ equivariant,
$$
s^\ast \omega = Ad(s) \omega
$$
for all $s\in S$.

Set $W = V\oplus \R$. Filter this by
$$
W = W^0 \supset W^1 \supset W^2 =0
$$
where $W^1 = V$. Then $V\subset \End W$. It follows from (\ref{converse})
that $\omega$ defines a connection on $P\times V$ which is flat along the
leaves of the foliation $\F$ and descends to a flat bundle over $X$. The
monodromy representation of this bundle is a homomorphism
$$
\tau:\pi_1(X,x_o) \to S\ltimes V \subset \Aut W.
$$
It follows from the monodromy formula (\ref{transp}) that $\tau$ takes the
class of the loop $\gamma$ to
$$
\left(\rho(\gamma),\int_\gammatilde w\right) \in S\ltimes V.
$$
The result follows.
\end{proof}

\section{Naturality}

Suppose that $\pi_X$ and $\pi_Y$ are two groups, and that
$\rho_X : \pi_X \to S_X$ and $\rho_Y : \pi_Y \to S_Y$ are homomorphisms
into the $F$-points of reductive algebraic groups, each with Zariski
dense image. We have the two corresponding relative completions
$$
\rhotilde_X : \pi_X \to \G_X \text{ and } \rhotilde_Y : \pi_Y \to \G_Y.
$$
Fix an algebraic group homomorphism $\Psi : S_X \to S_Y$.

\begin{proposition}\label{induced}
If $\psi : \pi_X \to \pi_Y$ is a homomorphism such that the diagram
$$
\begin{CD}
\pi_X @>{\rho_X}>> S_X \cr
@V{\psi}VV    @VV{\Psi}V \cr
\pi_Y @>{\rho_Y}>> S_Y
\end{CD}
$$
commutes, then there is a canonical homomorphism
$\Psihat : \G_X \to \G_Y$ such that the diagram
$$
\begin{CD}
\pi_X @>{\rhotilde_X}>> \G_X \cr
@V{\psi}VV    @VV{\Psihat}V \cr
\pi_Y @>{\rhotilde_Y}>> \G_Y
\end{CD}
$$
\end{proposition}

\begin{proof}
Let $\Psi^\ast \G_Y$ be the pullback of $G_Y$ along $\Psi$:
$$
\begin{CD}
\Psi^\ast \G_Y @>>> S_X \cr
@VVV @VV{\Psi}V \cr
\G_Y @>>> S_Y.
\end{CD}
$$
This group is an extension of $S_X$ by the prounipotent radical
of $\G_Y$.
The homomorphisms $\pi_X \to S_X$ and $\pi_X \to \pi_Y \to \G_Y$
induce a homomorphism $\pi_X \to \Psi^\ast \G_Y$. By the universal
mapping property of $\rhotilde_X : \pi_X \to \G_X$, there is a
homomorphism $\G_X \to \Psi^\ast\G_Y$ which extends the homomorphism
$\pi_X \to \Psi^\ast\G_Y$. The sought after homomorphism $\Psihat$
is the composite $\G_X \to \Psi^\ast\G_Y \to \G_Y$.
\end{proof}

Next we explain how to realize $\Psihat$ using the bar construction.
Suppose that $(X,x)$ and $(Y,y)$ are two pointed manifolds. Denote
$\pi_1(X,x)$ and $\pi_1(Y,y)$ by $\pi_X$ and $\pi_Y$, respectively.
Suppose that $f:(X,x) \to (Y,y)$ is a smooth map which induces the
homomorphism $\psi : \pi_X \to \pi_Y$ on fundamental groups. Denote
the principal bundles associated to $\rho_X$ and $\rho_Y$ by $P_X \to X$
and $P_Y \to Y$. We have the d.g.a.s
$$
\Efin^\dt(X,\O(P_X)) \text{ and } \Efin^\dt(Y,\O(P_Y)).
$$
Since the diagram in Proposition~\ref{induced} commutes, $f$ and $\psi$
induce a d.g.a.\ homomorphism
$$
(f,\phi)^\ast : \Efin^\dt(Y,\O(P_Y)) \to \Efin^\dt(X,\O(P_X)).
$$
This homomorphism respects the augmentations
induced by $x \in X$ and $y \in Y$, and therefore induces
a d.g.~Hopf algebra homomorphism
$$
B(\Efin^\dt(Y,\O(P_Y))_{\ytilde,(y)}) \to
B(\Efin^\dt(X,\O(P_X))_{\xtilde,(x)}).
$$
This induces a homomorphism
\begin{equation}\label{induced_homom}
\Gdr(X,x) \to \Gdr(Y,y)
\end{equation}
after taking $H^0$ and then $\spec$.

\begin{proposition}
Under the canonical identifications of $\Gdr(X,x)$ with $\G(X,x)$ and
$\Gdr(Y,y)$ with $\G(Y,y)$, the homomorphism
$\Psihat : \G(X,x) \to \G(Y,y)$
corresponds to the homomorphism
(\ref{induced_homom}).
\end{proposition}

\begin{proof}
If $\gamma$ is a loop in $X$ based at $x$,
$w_1,\dots,w_r \in \Efin^\dt(Y,\O(P_Y))$, and $\phi \in \O(S_Y)$, then
$$
\int_{f\circ \gamma}\left(w_1 \dots w_r | \phi\right)
\int_{\gamma}\left((f,\phi)^\ast w_1 \dots (f,\phi)^\ast w_r |
\Psi^\ast\phi\right).
$$
It follows that (\ref{induced_homom}) is the homomorphism $\Psihat$
induced by $f$ and $\psi$.
\end{proof}

\section{Relative Completion of the Fundamental Groupoid}

In this section we explain how the fundamental groupoid of $X$
can be completed with respect to $\rho : \pi_1(X,x) \to S$ and we
give a de~Rham construction of it. In the unipotent case, the de~Rham
theorem is implicit in Chen's work \cite{chen}, and is described
explicitly in \cite{hain-zucker}.

Recall that the fundamental groupoid $\pi(X)$ of a topological space $X$
is the category whose objects are the points of $X$ and whose morphisms
from $a \in X$ to $b \in X$ are homotopy classes $\pi(X;a,b)$ of paths
$[0,1] \to X$ from $a$ to $b$. We can think of $\pi(X)$ as a torsor over
$X\times X$; the fiber over $(a,b)$ being $\pi(X;a,b)$. Observe that
there is a canonical isomorphism between the  fiber over $(a,a)$ and
$\pi_1(X,a)$. The torsor is the one over $X\times X$ corresponding to
the representation
\begin{equation}\label{action}
\pi_1(X\times X,(a,a)) \cong \pi_1(X,a)\times \pi_1(X,a) \to \Aut \pi_1(X,a)
\end{equation}
where
$$
(\gamma, \mu) \to \left\{g \mapsto \gamma^{-1}g \mu\right\}.
$$

As in previous sections, $X$ will be a connected smooth manifold and
$x_o$ a distinguished base point. Suppose, as before, that
$\rho : \pi_1(X,x_o) \to S$ is a Zariski dense homomorphism to a reductive
real algebraic group. Denote the completion of $\pi_1(X,x_o)$  relative to
$\rho$ by $\pi_1(X,x_o) \to \G$. The representation (\ref{action})
extends to a representation
$$
\pi_1(X\times X,(x_o,x_o)) \cong \pi_1(X,x_o)\times \pi_1(X,x_o) \to \Aut \G
$$
Denote the corresponding torsor over $X\times X$ by $\bG$. This is easily
seen to be a torsor of real proalgebraic varieties. Denote the fiber
of $\bG$ over $(a,b)$ by $\G_{a,b}$. There is a canonical map
$$
\pi(X;a,b) \to \G_{a,b}
$$
which induces a map of torsors. It follows from standard arguments that,
for all $a$, $b$ and $c$ in $X$, there is a morphism of proalgebraic
varieties
\begin{equation}\label{mult_map}
\G_{a,b} \times \G_{b,c} \to \G_{a,c}
\end{equation}
which is compatible with the multiplication map
$$
\pi(X;a,b) \times \pi(X;b,c) \to \pi(X;a,c).
$$
An efficient way to summarize the properties of $\bG$ and the
multiplication maps is to say that they form a category (in fact, a
groupoid) whose objects are the elements of $X$ and where $\Hom(a,b)$
is $\G_{a,b}$ with composition defined by (\ref{mult_map}). In addition,
the natural map $\pi(X;a,b) \rightsquigarrow \G_{a,b}$ from the fundamental
groupoid of $X$ to this category is a functor. We shall call this functor
the {\it relative completion of the fundamental groupoid of $X$ with respect
to $\rho$}.\footnote{We shall see that the torsor $\bG$ is independent
of the choice of base point $x_o$, so it may have been better to call
$\bG$ the completion of the fundamental groupoid relative to the principal
bundle $P$.}
Our goal is to give a description of it in terms of differential
forms.

We also have the torsor $\bP$ over $X\times X$ associated to the
representation
$$
\pi_1(X\times X,x_o) \cong \pi_1(X,x_o)\times \pi_1(X,x_o) \to \Aut S
$$
where
$$
(\gamma, \mu)\to\left\{g \mapsto\rho(\gamma)^{-1}g\rho(\mu)\right\}.
$$
given by $\rho$. Denote the fiber of $\bP$ over $(a,b)$ by $\cP_{a,b}$.
As above, we have a category whose objects are the points of $X$ and
where $\Hom(a,b)$ is $\cP_{a,b}$. There is also a functor from the
fundamental groupoid of $X$ to this category which is the identity on
objects. Denote the restriction of this torsor to $\{a\}\times X$
by $\bP_{a,\omit}$. We shall view this as a torsor over $(X,a)$. The
image $\id_a$ of the identity in $\pi_1(X,a)$ in $\cP_{a,a}$ gives a
canonical lift of the base point $a$ of $X$ to $\bP_{a,\omit}$.
Observe that $\bP_{x_o,\omit}$ is the principal $S$ bundle $P$ used in
the construction of $\G$.

For each $a \in X$, we can from the corresponding local system
$\O(P_{a,\omit})$ whose fiber over $b\in X$ is the coordinate ring
of $\cP_{a,b}$. We can form the complex
$$
\Efin^\dt(X,\O(\cP_{a,\omit})) = \varinjlim E^\dt(X,\M)
$$
where $\M$ ranges over all finte dimensional sub-local systems of
$\O(\cP_{a,\omit})$. This has
augmentations
$$
\epsilon_{a} : \Efin^\dt(X,\O(\cP_{a,\omit})) \to \R
$$
and
$$
\delta_{a,b} : \Efin^\dt(X,\O(\cP_{a,\omit})) \to \O(\cP_{a,b})
$$
given by evaluation at $\id_a$ and on the fiber over $b$, respectively.
We view $\R$ as a right $\Efin^\dt(X,\O(\cP_{a,\omit}))$ module via
$\epsilon_a$ and $\O(\cP_{a,b})$ as  a left
$\Efin^\dt(X,\O(\cP_{a,\omit}))$ module via $\delta_{a,b}$.
We can therefore form the two sided bar construction
$$
B(\Efin^\dt(X,\O(\cP_{a,\omit}))_{\id_a,(b)}) :=
B(\R,\Efin^\dt(X,\O(\cP_{a,\omit})),\O(\cP_{a,b}))
$$
Define
\begin{multline}\label{comult}
B(\Efin^\dt(X,\O(\cP_{a,\omit}))_{\id_a,(c)}) \\ \to
B(\Efin^\dt(X,\O(\cP_{a,\omit}))_{\id_a,(b)}) \otimes
B(\Efin^\dt(X,\O(\cP_{b,\omit}))_{\id_b,(c)})
\end{multline}
by
\begin{multline*}
\Delta : [w_1|\dots|w_r]\phi \mapsto \cr
\sum_i [w_1|\dots |w_i]
\left(\sum \psi_i^{(k_i)}\dots\psi_r^{(k_r)}\phi'\right)
\otimes [w_{i+1}^{(k_i)}|\dots|w_r^{(k_r)}]\left(\sum \phi''\right)
\end{multline*}
where
$$
\Delta_S(\phi) = \sum \phi'\otimes \phi''
$$
is the map $\O(P_{a,c})\to \O(P_{a,b})\otimes\O(P_{b,c})$ dual to
the multiplication map $P_{a,b}\times P_{b,c}\to P_{a,c}$; and the map
$$
\Efin^\dt(X,\O(P_{a,\omit})) \to
\O(P_{a,b}) \otimes \Efin^1(X,\O(P_{b,\omit}))
$$
induced by multiplication $P_{a,b} \times \bP_{b,\omit} \to \bP_{a,\omit}$
takes $w_j$ to
$$
\sum \psi_j^{(k_j)} \otimes w_j^{(k_j)}.
$$

Definition~\ref{defn} generalizes:

\begin{definition}
For $\gamma$ a path in $X$ from $a$ to $b$, $\phi\in \O(P_{a,b})$ and
$w_1,\dots,w_r$ elements of $\Efin^1(X,\O(P_{a,\omit}))$, we define
$$
\int_\gamma \left(w_1 \dots w_r | \phi\right)
= \phi(\gammatilde(1))\int_{\gammatilde}w_1\dots w_r
$$
where $\gammatilde$ is the unique lift of $\gamma$ to a horizontal
section of $\bP_{a,\omit}$ which begins at $\id_a \in P_{a,a}$.
\end{definition}

There is an analogous extension of the definition of higher
iterated integrals (\ref{higherdef}) to this situation. As in
that case, one has a d.g.\ algebra homomorphism
$$
B(\Efin^\dt(X,\O(\cP_{a,\omit}))_{\id_a,(b)}) \to E^\dt(P_{a,b}X)
$$
to the de~Rham complex of $P_{a,b}X$, the space of paths in $X$
from $a$ to $b$. It is defined by
$$
[w_1|\dots|w_r]\phi \mapsto \int(w_1\dots w_r|\phi).
$$
By taking a homotopy class $\gamma \in \pi(X;a,b)$ to the ideal of
functions that vanish on it, we obtain a function
$$
\pi(X;a,b) \to
\spec H^0(B(\Efin^\dt(X,\O(\cP_{a,\omit}))_{\id_a,(b)})).
$$

\begin{theorem}\label{gpoid_dr}
This function gives a natural algebra isomorphism
$$
\O(\G_{a,b}) \cong
H^0(B(\Efin^\dt(X,\O(\cP_{a,\omit}))_{a,(b)})).
$$
Moreover, the map
$$
\O(\G_{a,c}) \to \O(\G_{a,b})\otimes \O(\G_{b,c})
$$
induced by (\ref{mult_map}) corresponds to (\ref{comult}) under this
isomorphism.
\end{theorem}

\begin{proof}[Sketch of Proof]
Define $\Gdr_{a,b}$ by
$$
\Gdr_{a,b} = \spec H^0(B(\Efin^\dt(X,\O(\cP_{a,\omit}))_{a,(b)})).
$$
The coproduct above induces morphisms
$$
\Gdr_{a,b} \times \Gdr_{b,c} \to \Gdr_{a,c}.
$$
We therefore have a groupoid whose objects are the points of $X$
and where $\Hom(a,b)$ is $\Gdr_{a,b}$ and a function
$$
\pi(X;a,b) \to \Gdr_{a,b}.
$$
This map is easily seen to be compatible with path multiplication
(use the generalization of the last property of (\ref{props})),
and therefore a functor of groupoids. Since $X$ is connected, it
suffices to prove that
$\Gdr_{a,b}$ is isomorphic to $\G_{a,b}$ for just one pair of points
$a,b$ of $X$. But these are isomorphic in the case $a=b=x_o$ by
Theorem~\ref{main}.
\end{proof}

\section{Hodge Theory}
\label{hodge_str}

Now suppose that $X$ is a smooth complex algebraic variety (or the
complement of a normal crossings divisor in a compact K\"ahler manifold)
and that $\V$ is an admissible variation of Hodge structure over $X$.
Denote the semisimple group associated to the fiber $V_o$ over
the base point $x_o\in X$ by $S$.  This is the ``orthogonal'' group
$$
S = \Aut(V_o,\langle\blank,\blank\rangle)
$$
associated to the polarization $\langle\blank,\blank\rangle$. It is
semi-simple. Suppose that the image of the monodromy representation
$$
\rho : \pi_1(X,x_o) \to S
$$
is Zariski dense. Denote the completion of $\pi_1(X,x_o)$ relative to
$\rho$ by
$$
\rhotilde : \pi_1(X,x_o) \to \G(X,x_o).
$$

\begin{theorem}\label{hodge}
Under these assumptions, the coordinate ring $\O(\G(X,x_o))$ of the
completion of $\pi_1(X,x_o)$ with respect to $\rho$ has a canonical
real mixed Hodge structure with weights $\ge 0$ for which the product,
coproduct, antipode and the natural inclusion
$$
\O(S) \hookrightarrow \O(\G(X,x_o))
$$
are all morphisms of mixed Hodge structure. Moreover the canonical
homomorphism $\G(X,x_o) \to S$ induces an isomorphism
$\Gr^W_0\O(\G(X,x_o)) \cong \O(S)$.
\end{theorem}

Denote the Lie algebra of $S$ by $\s$. This has a canonical Hodge
structure of weight 0.
The following result is an important corollary of the proof of Theorem~%
\ref{hodge}. It follows immediately from the theorem and the standard
description of the Lie algebra of an affine algebraic group given at the
end of Section~\ref{coord}.

\begin{corollary}\label{hodge-lie}
Under the assumptions of the theorem, the Lie algebra $\g(X,x_o)$ of
$\G(X,x_o)$ has a canonical MHS with weights $\le 0$, and the homomorphism
$\g(X,x_o) \to \s$ is a morphism of MHS which induces an isomorphism
$$
\Gr^W_0 \g \cong \s.
$$
In particular, there is a canonical
MHS  with weights $< 0$ on $\u(X,x_o)$, the Lie algebra of the prounipotent
radical of $\G(X,x_o)$. \qed
\end{corollary}

The principal assertion of Theorem~\ref{hodge} is a special case of the
following result when $a=b=c=x_o$.

\begin{theorem}\label{groupoid}
With $X$, $\V$ and $S$ as above, if $a,\, b\in X$, then the coordinate
ring $\O(\G_{a,b})$ of the completion of $\pi(X;a,b)$ relative to $\rho$
has a canonical mixed Hodge structure with weight $\ge 0$ and whose
multiplication is a morphism of MHS. If $a$, $b$ and $c$ are three points
of $X$, then the map
$$
\O(\G_{a,c}) \to \O(\G_{a,b}) \otimes \O(\G_{b,c})
$$
induced by path multiplication is a morphism of MHS. Moreover, the
mixed Hodge structure on $\O(\G_{a,b})$ depends only on the variation
$\V$ and not on the choice of the base point $x_o$.
\end{theorem}

Because of the last assertion, it may be more appropriate to say
that {\it $\G_{a,b}$ is the completion of $\pi(X;a,b)$ with respect to
the variation $\V$.}

The reader is assumed to be familiar with the basic methods for
constructing mixed Hodge structures on the cohomology of bar
constructions as described in \cite[\S3]{hain:dht}.
In the previous section we showed how to express $\O(\G_{a,b})$ as the
0th cohomology group of a suitable reduced bar construction. So in
order to show that it has a canonical MHS we need only find a suitable
augmented, multiplicative mixed Hodge complex $\A^\dt$ which is
quasi-isomorphic to $\Efin^\dt(X,\O(\cP_{a,\omit}))$. To do this, we
shall use the work of M.~Saito on Hodge modules.

First, some notation: Assume that $X = \Xbar - D$, where
$\Xbar$ is a compact K\"ahler manifold and $D$ is a normal
crossings divisor. Denote the inclusion $X \hookrightarrow \Xbar$ by
$j$. Denote Deligne's canonical extension of
$\V\otimes\O_X$ to $\Xbar$ by $\Vbar$. Saito proves that there is a
Hodge module over $\Xbar$ canonically associated to $\V$, whose
complex part is a bifiltered $D$-module $(M,W_\dt,F^\dt)$, and whose
real part is $Rj_\ast\V_\R$ endowed with a suitable weight filtration.
There is a canonical inclusion
$$
\Omega^\dt_\Xbar(\Xbar\log D)\otimes_\O \Vbar
\hookrightarrow M\otimes_\O \Omega_\Xbar^\dt.
$$
Saito defines Hodge and weight filtrations on
$\Omega^\dt_\Xbar(\Xbar\log D)\otimes_\O \Vbar$
by restricting those of $M$. The Hodge filtration is simply the
tensor product of those of $\Omega^\dt_\Xbar(\Xbar\log D)$
and $\Vbar$. The weight filtration is more difficult to describe.

\begin{theorem}[Saito \protect{\cite[(3.3)]{saito}}]\label{saito:mhc}
The pair
\begin{equation}\label{mhc}
M^\dt(X,\V) :=
((Rj_\ast\V_\R, W_\dt),
(\Omega^\dt_\Xbar(\Xbar\log D)\otimes_\O \Vbar, F^\dt, W_\dt))
\end{equation}
is a cohomological mixed Hodge complex whose cohomology is canonically
isomorphic to $H^\dt(X,\V)$. \qed
\end{theorem}

We can therefore obtain a mixed Hodge complex (MHC) which computes
$H^\dt(X,\V)$ by taking the standard fine resolution of these sheaves
by $C^\infty$ forms. (So the complex part of this will be the $C^\infty$
log complex $E^\dt(\Xbar\log D,\Vbar)$ with suitable Hodge and weight
filtrations.)

To apply Saito's machinery to the current situation, we will need
to know that $\O(\cP_{x_o,\omit})$ is a direct limit of admissible
variations over $X$.

\begin{lemma}\label{loc_sys}
The local system associated with an irreducible representation
of $S(\R)$ underlies an admissible variation of Hodge structure
over $X$. These structures are compatible with the decomposition
of tensor products. Moreover, these variations are unique up to Tate
twist.
\end{lemma}

\begin{proof}
The connected component of the identity of $S(\R)$ is a real form
of $Sp_n(\C)$ when $\V$ has odd
weight, and $SO_n(\C)$ when $\V$ has even weight. In both cases
each irreducible representations of the complex group can be
constructed by applying a suitable Schur functor the the
fundamental representation and then taking the intersection of the
kernels of all contractions with the polarization.
(This is Weyl's construction of the irreducible representations; it
is explained, for example, in \cite[\S17.3,\S19.2]{fulton-harris}.)
Since these operations preserve variations of Hodge structure, it
follows that a local system corresponding to an irreducible
representation of $S$ underlies a variation of Hodge structure.
Since the monodromy representation of $\V$ is Zariski dense, the
structure of a polarized variation of Hodge structure on this local
system is unique up to Tate twist. (Cf.\ the proof of
\cite[(9.1)]{hain:normal}.)
\end{proof}

This, combined with (\ref{decomp}) yields:

\begin{corollary}
With our assumptions, $\O(\cP_{x_o,\omit})$ is a direct limit of
admissible variations of Hodge structure over $X$ of weight 0,
and the multiplication map is a morphism. \qed
\end{corollary}

\begin{corollary}
For each $b\in X$, there is a canonical Hodge structure on
$\O(P_{x_o,b})$. \qed
\end{corollary}

Combining (\ref{loc_sys}) with (\ref{h1_comp}), we obtain:

\begin{corollary}\label{varmhs}
The local system over $X$ whose fiber over $x\in X$ is
$H^1(\U(X,x))$ is an admissible variation of mixed Hodge structure
whose weights are positive.
\end{corollary}

Using Saito's machine \cite{saito}, we see that there is a MHC $\A^\dt$
which is quasi-isomorphic to $\Efin^\dt(X,\O(\cP_{a,\omit}))$. The
complex part of this MHC is simply the complex of $C^\infty$ forms
on $\Xbar$ with logarithmic singularities along $D$ and which have
coefficients in the canonical extension $\Obar$ of $\O(\cP_{a,\omit})$
to $\Xbar$.
The Hodge filtration is the obvious one inducedd by the Hodge filtration
of $\Obar$ and the Hodge filtration of forms on $\Xbar$. The weight
filtration is not so easily described, and we refer to Saito's paper
for that.

We need to know that the multiplication is compatible with the
Hodge and weight filtrations. This follows from the next
result.

\begin{proposition}
If $\V_1\otimes\V_2 \to \W$ is a pairing of admissible
variations of Hodge structure over $X$, then the multiplication
map
$$
M^\dt(X,\V_1) \otimes M^\dt(X,\V_2) \to M^\dt(X,\W)
$$
is a morphism of cohomological mixed Hodge complexes.
\end{proposition}

\begin{proof}
This follows immediately from the naturality of Saito's construction,
its compatibility with exterior tensor products,
and the fact that admissible variations of Hodge structure are closed
under exterior products --- use restriction to the diagonal.
\end{proof}

There are two augmentations
$$
\R \leftarrow \A^\dt \to \O(P_{a,b})
$$
corresponding to the inclusions $P_{a,b}\hookrightarrow \cP_{a,\omit}$
and $\id_a \in P_{a,a}$, and these are compatible with
all filtrations. It follows from \cite[(3.2.1)]{hain:dht}, (\ref{qism})
and (\ref{gpoid_dr}) that
$$
B(\R,\A^\dt,\O(P_{a,b}))
$$
is a MHC whose $H^0$ is isomorphic to $\O(\G_{a,b})$. Moreover, the
multiplication is compatible
with the Hodge and weight filtrations. Consequently,
$$
\O(\G_{a,b}) \cong H^0(B(\R,\A^\dt,\O(P_{a,b})))
$$
has a canonical MHS and its multiplication is a morphism
of MHS. Since the MHS on $P_{a,b}$ depends only on $\V$ and
not on $x_o$, the same is true of the MHS on $\O(\G_{a,b})$.

If $a$, $b$, and $c$ are three points of $X$, then it follows directly
from the definitions that the map
$$
B(\R,\A^\dt,\O(P_{a,c})) \to
B(\R,\A^\dt,\O(P_{a,b}))\otimes B(\R,\A^\dt,\O(P_{b,c}))
$$
corresponding to path multiplication is a morphism of MHCs. It
follows that the induced map
$$
\O(\G_{a,c}) \to \O(\G_{a,b}) \otimes \O(\G_{b,c})
$$
is a morphism of MHS. This completes the proof of Theorem~\ref{groupoid};
Theorem~\ref{hodge} follows by taking $a=b=x_o$ except for the assertion
that $\O(S) \hookrightarrow \O(\G)$ is a morphism of MHS. This follows
as this is induced by the natural inclusion
$$
\O(S) \hookrightarrow B(\R,\A^\dt,\O(S)),
$$
that takes $\phi$ to $[\blank]\phi$. It is a morphism of MHCs. This
completes the proofs of Theorems~\ref{hodge} and \ref{groupoid}. \qed

We now turn our attention to the variation of the Hodge filtration.
Suppose that $X$ and $S$ are as above. Consider the real local
system over $X\times X$ whose fiber over $(a,b)$ is $\O(\G_{a,b})$.
Denote it by $\bO$.
Next we establish that this underlies a ``pre-variation of MHS.''
Denote by $\F^p\bO$the subset of the associated complex local system
with fiber $F^p\O(\G_{a,b}(\C))$ over $(a,b)$. Denote by $W_m\bO$
the subset of $\bO$ with fiber $W_m\O(\G_{a,b})$ over $b$.

\begin{theorem}\label{pre_var}
The subset $W_m\bO$ is a flat sub-bundle of $\bO$, and $F^p\bO$ is a
holomorphic sub-bundle of $\bO_\C$ whose corresponding sheaf of
sections $\F^p$ satisfies Griffiths transversality:
$$
\nabla : \F^p \to \F^{p-1}\otimes \Omega^1_X.
$$
\end{theorem}

\begin{proof}[Sketch of Proof]
We will prove the result for the restriction $\bO_a$ of $\bO$ to
$\{a\}\times X$. The result for the restriction of $\bO$ to
$X\times \{b\}$ is proved similarly. The general result then
follows as the tangent spaces of $\{a\}\times X$ and $X\times \{b\}$
span the tangent space of $X\times X$ at $(a,b)$.

First we need a formula for the connection on $\bO_a$ at the point
$b$ in $X$. Fix a path $\gamma$ in $X$ from $a$ to $b$.
Suppose that $\mu : [-\epsilon, \epsilon] \to X$ is a smooth path
with $\mu(0) = b$. For $s\in [-\epsilon, \epsilon]$ let
$\gamma_s : [0,1] \to X$ be the piecewise smooth path obtained by
following $\gamma$ and then $\mu$ from $t=0$ to $t=s$. Suppose that
$w_1,\dots, w_r$ are in $\Efin^\dt(X,\O(P_{a,b}))$. Suppose that
$U$ is a contractible neighbourhood of $b$ in $X$. With respect to
a flat trivialization of the restriction of $\cP_{a,\omit}$ to $U$,
we have
$$
w_r|_U = \sum_j w_r^j\otimes \psi_j
$$
where $w_r^j\in E^1(U)$ and $\psi_j \in \O(P_{a,b})$.

It follows from the analogue of (\ref{props}) in this situation that
$$
\frac{d}{ds}\bigg\vert_{s=0} \int_{\gamma_s}(w_1\dots w_r|\phi)
= \sum_j \int_\gamma (w_1\dots w_{r-1}|\phi \psi_j)
\langle w_r^j,\dot{\mu}(0)\rangle.
$$
The restriction of the connection on $\bO_a$ to the stalk at $b$ is
therefore induced by the map
$$
[w_1|\dots |w_r]\phi \mapsto
\sum_j [w_1|\dots |w_{r-1}]\phi \psi_j \otimes w_r^j
$$
on the bar construction.
The flatness of the weight filtration follows immediately from
the definition of the weight filtration on the bar construction.
Further, if $(z_1,\dots, z_n)$ is a holomorphic coordinate in $X$
centered at $b$, then it follows immediately from the definition of
the Hodge and weight filtrations on $\O(\G_{a,b}(\C))$ and the formula
for the connection that
$$
\nabla_{\partial/\partial \zbar_k} : \F^p \to \F^p
$$
for each $k$, so that the Hodge filtration varies holomorphically
at $b$. Similarly, on the stalk of $F^p$ at $b$ we have
$$
\nabla_{\partial/\partial z_k} : \F^p \to \F^{p-1}
$$
as each $w_r^j$ contributes at most 1 to the Hodge filtration
of $\O(\G_{a,b})$.
\end{proof}

When $X$ is compact, we have:

\begin{corollary}\label{good_varn}
If $X$ is a compact K\"ahler manifold, then $\bO$ is an admissible
variation of MHS over $X\times X$. \qed
\end{corollary}

In order to prove the corresponding result in the non-compact case,
it is necessary to study the asymptotic behaviour of $\bO$. I plan
to consider this in a future paper.

We now consider naturality. Suppose that $X$ and $\V$ are as
above, and that $Y$ is a smooth variety and that $\W$ is an
admissible variations of Hodge structure over $Y$. We will now
denote $S$ by $S_X$:
$$
S_X = \Aut(V_o,\langle\blank,\blank\rangle).
$$
Set
$$
S_Y = \Aut(W_o,\langle\blank,\blank\rangle)
$$
where $W_o$ denotes the fiber of $\W$ over $y_o$. Suppose that
the monodromy representation
$$
\rho_Y : \pi_1(Y,y_o) \to S_Y
$$
has Zariski dense image. Denote the completion of $\pi_1(X,x_o)$
with respect to $\rho_X:\pi_1(X,x_o) \to S_X$ by $\pi_1(X,x_o)\to \G_X$,
and the completion of $\pi_1(Y,y_o)$ with respect to $\rho_Y$ by
$\pi_1(Y,y_o) \to \G_Y$.

Suppose that $f:(Y,y_o) \to (X,x_o)$ is a morphism of varieties, and
that we have fixed an inclusion
$$
\End \V \hookrightarrow \End f^\ast \W
$$
of variations of Hodge structure. This fixes a group homomorphism
$$
\Psi : S_X \hookrightarrow S_Y
$$
that is compatible with the Hodge theory.
By (\ref{induced}), there is a homomorphism $\Psihat : \G_X \to \G_Y$
such that the diagram
$$
\begin{CD}
\pi_1(X,x_o) @>>> \G_X @>>> S_X \cr
@V{f_\ast}VV @VV{\Psihat}V @VV{\Psi}V \cr
\pi_1(Y,y_o) @>>> \G_Y @>>> S_Y \cr
\end{CD}
$$
commutes.

\begin{theorem}\label{hodge_nat}
Under these hypotheses, the induced map
$$
\Psihat^\ast : \O(\G_Y) \to \O(\G_X)
$$
is a morphism of MHS.
\end{theorem}

\begin{proof}
First, choose smooth compactifications $\Xbar$ of $X$ and $\Ybar$
of $Y$ such that $X=\Xbar - D$ and $Y = \Ybar - E$, where $D$ and
$E$ are normal crossings divisors. We may choose these such that $f$
extends to a morphism $\Xbar \to \Ybar$, which we shall also denote
by $f$.

Denote by $P_X$ the flat left $S_X$ principal bundle over $X$
associated the the
representation of $\pi_1(X,x_o)$ on $S_X$ via $\rho_X$. Denote the
analogous principal $S_Y$ principal bundle over $Y$ by $P_Y$.
Associated to these we have the local systems $\O(P_X)$ over $X$
and $\O(P_Y)$ over $Y$.

The construction above gives multiplicative mixed Hodge complexes
$$
\A^\dt(X,\O(P_X)),\quad \A^\dt(X,f^\ast\O(P_Y)),
\text{ and } \A^\dt(Y,\O(P_Y))
$$
which compute the canonical mixed Hodge structures on
$$
H^\dt(X,\O(P_X)),\quad H^\dt(X,f^\ast\O(P_Y)),\text{ and }
H^\dt(X,f^\ast\O(P_Y))
$$
respectively. The map $f$ induces a morphism
$$
\A^\dt(Y,\O(P_Y)) \to \A^\dt(X,f^\ast\O(P_Y))
$$
of MHCs, while the inclusion $S_X \hookrightarrow S_Y$ induces a
morphism
$$
\A^\dt(X,f^\ast\O(P_Y)) \to \A^\dt(X,\O(P_X))
$$
of MHCs. The composition of these corresponds to the induced map
$$
\Efin^\dt(Y,\O(P_Y)) \to \Efin^\dt(X,\O(P_X))
$$
induced by $f$ under the canonical quasi-isomorphisms. It follows that
the induced map
$$
B(\R,\A^\dt(Y,\O(P_Y)),\O(P_Y)) \to B(\R,\A^\dt(Y,\O(P_X)),\O(P_X))
$$
is a morphism of MHCs and that the induced the map
$$
f^\ast : \O(\G_Y) \to \O(\G_X)
$$
on $H^0$ is the ring homomorphism induced by $f$. The result follows.
\end{proof}

\begin{remark}\label{extended}
Suppose that $\V$ is an admissible variation of Hodge structure over
the complement $X$ of a normal crossings divisor in a compact K\"ahler
manifold. We will say that the pair $(X,\V)$ is {\it neat} if the Zariski
closure $S$ of the image of the monodromy map
$$
\rho : \pi_1(X,x_o) \to \Aut(V_o,\langle\blank,\blank\rangle)
$$
is semi-simple, and that the canonical MHS on the coordinate ring of
$$
\Aut(V_o,\langle\blank,\blank\rangle)
$$
induces one on $S$.  For example, every variation where $S$ is finite
is neat. I believe that every admissible $(X,\V)$ is neat, but have not
yet found a proof.

The results (\ref{hodge}), (\ref{hodge-lie}), (\ref{groupoid}),
(\ref{varmhs}), (\ref{pre_var}), (\ref{good_varn}), (\ref{hodge_nat})
and their proofs are valid with the assumption that $\im \rho$ be
Zariski dense in $\Aut(V_o,\langle\blank,\blank\rangle)$ replaced
by the assumption that the pairs $(X,\V)$ and $(Y,\W)$ be neat.
\end{remark}

The following is an application suggested by Ludmil Kartzarkov.

\begin{theorem}
Suppose that $X$ is a compact K\"ahler manifold and that $\V$ is
a polarized variation of Hodge structure over $X$ whose monodromy
representation $\rho$ has Zariski dense image. Then the prounipotent
radical of the completion of $\pi_1(X,x_o)$ relative to $\rho$ has
a quadratic presentation.
\end{theorem}

\begin{proof}
It is well known that if $X$ is compact K\"ahler and $\V$ is a
polarized variation of Hodge structure over $X$ of weight $m$,
then $H^k(X,\V)$ has a pure Hodge structure of weight $k+m$.
In particular, as the variation $\O(P)$ over $X$ has weight
zero, $H^k(X,\O(P))$ is pure of weight weight $k$ for all $k$.

Denote the Lie algebra of the prounipotent radical of the relative
completion of $\pi_1(X,x_o)$ by $\u$.
It follows from (\ref{def}) and (\ref{main}) that $\u$ is the Lie
algebra canonically associated to the d.g.a.\ $\Efin^\dt(X,\O(P))$
by rational homotopy theory (either via Sullivan's theory of minimal
models, or via Chen's theory as the dual of the indecomposables
of the bar construction on $\Efin^\dt(X,\O(P))$.)  There are canonical
maps
$$
H^1(\u) \cong H^1(X,\O(P))\text{ and }
H^2(\u) \hookrightarrow H^2(X,\O(P)).
$$
These are morphisms of MHS \cite[(7.2)]{hain:torelli}. It follows
that $H^1(\u)$ is pure of weight 1 and $H^2(\u)$ is pure of weight 2.
The result now follows from \cite[(5.2),(5.7)]{hain:torelli}.
\end{proof}

\begin{remark}
It is not necessarily true that $\u$ is a quotient of the unipotent
completion of $\ker \rho$. A criterion for surjectivity is given in
\cite[(4.6)]{hain:comp}. For this reason it may not be easy to apply this
result in general situations without artificially restrictive hypotheses.
\end{remark}

\section{A Canonical Connection}
\label{connection}

For the time being, let $X$, $x_o$, $\V$, $\rho$, etc.\ be as in the
previous section. However, all groups and Lie algebras in this section
will be complex, and $\G$, $\U$, $\u$, etc.\
denote the {\em complex} points of the relative completion of
$\pi_1(X,x_o)$, its prounipotent radical, its Lie algebra, etc.
Denote the image of $\rho$ by $\Gamma$,
and the Galois covering of $X$ with Galois group $\Gamma$ by $X'$.
In this section, we show how the Hodge theory of $\G$ gives a
canonical (given the choice of $x_o$), $\Gamma$ invariant integrable
1-form
$$
\omega \in E^1(X')\comptensor \Gr^W_\dt\u
$$
on $X'$ which can be integrated to the canonical representation
$$
\rhotilde : \pi_1(X,x_o) \to S \ltimes \U \cong \G.
$$
Here $\comptensor$ denotes the completed tensor product, which
is defined below.

At the end of the section, we shall explain what this means when
$X$ is the complement of the discriminant locus in $\C^n$ and $S$ is
the symmetric group $\Sigma_n$. In this case, $X'$ is the
complement of the hyperplanes $x_i = x_j$ in $\C^n$, $\pi_1(X,x_o)$
is the classical braid group $B_n$, and
the form is
$$
\omega = \sum_{i<j} d\log(x_i-x_j).
$$

First, we shall define the completed tensor product $\comptensor$.
Suppose that $\u$ is a topological Lie algebra and that
$$
\u = \u^1 \supseteq \u^2 \supseteq \u^3 \supseteq \cdots
$$
is a base of neighbourhoods of 0. Suppose that
$E$ is a vector space. Define
$$
E\comptensor \u = \lim_{\leftarrow} E\otimes \u/\u^m.
$$
We can regard a graded Lie algebra $\u = \oplus_{m < 0} \u_m$ as a
topological Lie algebra where the basic neighbourhoods of 0 are
$$
\bigoplus_{l\le m} \u_l, \quad m < 0.
$$
The definition of completed tensor product therefore extends to
the case where $\u$ is graded. Finally, if $\u$ is a Lie algebra
in the category of mixed Hodge structures where $\u = W_{-1}\u$
which is complete with respect to the weight filtration, and if $E$
is a complex vector space, then there is a canonical isomorphism
$$
E\comptensor \u_\C \cong E \comptensor \Gr^W_\dt \u_\C
$$
as $\u_\C$ is canonically isomorphic to
$\prod \Gr^W_m\u_\C$. (Cf.\ \cite[(5.2)]{hain:torelli}.)

We view a principal bundles with structure group a
proalgebraic group to be the inverse limit of the principal
bundles whose structure groups are the finite dimensional
quotients of the proalgebraic group. A connection on a
principal bundle with proalgebraic structure group is
the inverse limit of compatible connections on the
corresponding bundles with finite dimensional structure
group.

\subsection{The unipotent case.}
We begin with the unipotent case, $S=1$. Suppose that $X$ is a
smooth manifold with distinguished base point $x_o$. Denote the complex
form of the unipotent completion of $\pi_1(X,x_o)$ by $\G$ and the
complex points of the completion of $\pi(X;x_o,x)$ by $\G_{x_o,x}$. The
family
$$
\left(\G_{x_o,x}\right)_{x\in X}
$$
forms a flat principal left $\G$ bundle over $X$ that we shall
denote by $\G_{x_o,\omit}$. Since the structure group is contractible,
(more precisely, an inverse limit of contractible groups), this bundle
has a section. Pulling back the canonical connection form, we obtain
an integrable connection form
$$
\omega \in E^\dt(X)\comptensor \g
$$
where $\g$ denotes the Lie algebra of $\G$. The monodromy
representation of this form is the monodromy of $\G_{x_o,\omit}$,
which is the canonical homomorphism
$$
\pi_1(X,x_o) \to \G.
$$

When $X$ is an algebraic manifold, there is a canonical choice of
section and therefore a canonical connection form. To see this,
note that for each $a \in X$, the weights on $\O(\G_{x_o,a})$ are
$\ge 0$ and that
$$
\Gr^W_0 \O(\G_{x_o,a}) \cong \C.
$$
Since there is a canonical ring isomorphism
$$
\O(\G_{x_o,a}) \cong
\bigoplus_{l\ge 0} \Gr^W_l \O(\G_{x_o,a})
$$
there is a canonical augmentation
$$
\O(\G_{x_o,a}) \to \C
$$
whose kernel is
$$
\bigoplus_{l > 0} \Gr^W_l \O(\G_{x_o,a}).
$$
This determines a canonical point in $\G_{x_o,a}$.
Since the family $\left\{\O(\G_{x_o,a})\right\}_{a\in X}$ is a variation
of MHS over $X$ (see \cite{hain-zucker}), these distinguished points
vary smoothly as $a$ varies. They therefore determine a smooth section
of $\G_{x_o,\omit}$. We therefore have a canonical integrable
1-form
$$
\omega \in E^1(X)\comptensor\g \cong E^1(X)\comptensor\Gr^W_\dt\g.
$$

\subsection{The general case}
The first step in doing this in general is to explain the necessary
constructions in the $C^\infty$ case. So suppose for the time being that
$X$ is a smooth manifold; $\rho$, $S$, $\G$ and $P\to X$ are as before.
We also have the torsor
$$
\G_{x_o,\omit} \to X.
$$
It is a flat principal left $\G$ bundle over $X$. There is a map
$$
\begin{CD}
\G_{x_o,\omit} @>\pi>> P \cr
@VVV					@VVV \cr
X @= X
\end{CD}
$$
of flat bundles. It is compatible with the canonical homomorphism
$\G \to S$. Denote by $X'$ the leaf of $P$ containing the distinguished
lift $\xtilde_o$ to $P$ of $x_o$. The projection $P\to X$ induces a
covering map $X' \to X$. It is the Galois covering corresponding to
$\ker \rho$. Define $\U_{x_o,\omit}$ to be the subset $\pi^{-1}X'$
of $\G_{x_o,\omit}$.
There is a natural projection $\U_{x_o,\omit}\to X'$ indued by $\pi$.
Note that the fiber of this over $\xtilde_o$ is $\U$, the prounipotent
radical of $\G$. Denote the fiber of $\U_{x_o,\omit}$ over $a\in X'$ by
$\U_{x_o,a}$.

Each point $a$ of $P$ determines an augmentation
$$
\epsilon_a : \Efin^\dt(X,\O(P)) \to \C.
$$
Given two points $a$ and $b$ of $P$, we may form the two sided
bar construction
\begin{equation}\label{bar}
B(\C,\Efin^\dt(X,\O(P)),\C)
\end{equation}
where the left hand $\C$ is viewed as a module over $\Efin^\dt(X,\O(P))$
via $\epsilon_a$, and the right hand $\C$ via $\epsilon_b$. We shall
denote the d.g.a.\ (ref{bar}) by $B(\Efin^\dt(X,\O(P))_{a,b})$

\begin{proposition}
Each $\U_{x_o,a}$ is a proalgebraic variety with coordinate ring
$$
\O(\U_{x_o,a}) \cong H^0(B(\Efin^\dt(X,\O(P))_{\xtilde_o,a})).
$$
Moreover, $\U_{x_o,\omit}\to X'$ is a principal $\U$ bundle with respect
to the natural $\U$ action on $\G_{x_o,\omit}$.
\end{proposition}

Choose a splitting $S\to \G$ of the natural homomorphism $\G \to S$.
This induces an isomorphism $\G\cong S\ltimes \U$. The splitting enables
us to lift the action of $S$ to $\G_{x_o,\omit}$ in such a way that
the projection $\G_{x_o,\omit}\to P$ is $S$ equivariant. Since $\Gamma$
is a subgroup of $S$, and since it preserves $X'\subset P$, it follows
that there is a natural left action of $\Gamma$ on $\U_{x_o,\omit}$ and
that, with respect to this action, the projection $\U_{x_o,\omit}$
is $\Gamma$ equivariant.

Denote the pullback of the extension
$$
1 \to \U \to \G \to S \to 1
$$
along $\Gamma \hookrightarrow S$ by $\G_\Gamma$. This is an extension
$$
1 \to \U \to \G_\Gamma \to \Gamma \to 1.
$$
The splitting $S\to \G$ induces a splitting $\Gamma \to \G_\Gamma$,
and therefore a semi-direct productu decomposition
$\G_\Gamma \cong \Gamma \ltimes \U$.
The image of the canonical homomorphism $\pi_1(X,x_o) \to \G$ lies in
$\G_\Gamma$. The composite $\U_{x_o,\omit} \to X' \to X$ is a
flat principal left $\G_\Gamma$ bundle over $X$. The associated monodromy
representation is the canonical homomorphism $\pi_1(X,x_o) \to \G_\Gamma$.
The monodromy therefore induces the canonical homomorphism
$$
\pi_1(X,x_o) \to \G \cong S \ltimes \U.
$$
Next, we explain that the pullback of this bundle to $X'$ is trivial,
and therefore given by an integrable 1-form.

\begin{proposition}
There is a $\Gamma$ equivariant section of $\U_{x_o,\omit}\to X'$.
\end{proposition}

\begin{proof}
The action of $\Gamma$ on $X'$ is free. It follows that the action
of $\Gamma$ on $\U_{x_o,\omit}$ is also free. Consequently, the
square
$$
\begin{CD}
\U_{x_o,\omit} @>>> \Gamma\backslash \U_{x_o,\omit} \cr
@VVV									@VVV \cr
X' @>>p> X
\end{CD}
$$
is a pullback square. Since the fibers of
$\Gamma\backslash \U_{x_o,\omit} \to X$
are connected, it has a $C^\infty$ section.  This section pulls back
to a $\Gamma$ invariant section of $\U_{x_o,\omit} \to X'$.
\end{proof}

Let $\Gamma$ act on $\U$ on the left via the adjoint action:
$$
Ad(\gamma) : u \mapsto \gamma u \gamma^{-1}.
$$
Then $\Gamma$ acts on $X'\times \U$ on the left via the diagonal
action. It follows from the previous result that the flat principal
bundle $\U_{x_o,\omit}\to X'$ has a $\Gamma$ invariant trivialization.
We therefore have a connection form
$$
\omega \in E^1(X')\comptensor \u.
$$

\begin{proposition}
This connection form satisfies $\gamma^\ast \omega = Ad(\gamma)\omega$
for all $\gamma \in \Gamma$.
\end{proposition}

\begin{proof}
Since the bundle is trivial, its (locally defined) sections can be
identified with (locally defined) functions $X' \to \U$. Since $\Gamma$
preserves the connection, we see that for each $\gamma\in \Gamma$ the
local section $u$ is flat if and only if the local section
$(\gamma^{-1})^\ast Ad(\gamma)(u)$ is flat. That is, $Ad(\gamma)(u)$
is flat if and only if $\gamma^\ast u$ is flat. The result now follows
from a standard and straight forward computation.
\end{proof}

Now suppose that $X$ is an algebraic manifold. We have to show
that this construction can be made canonical. Note that, given the
choice of the base point $x_o$, the only
choices made in the construction of $\omega$ were the choice of a
splitting of the homomorphism $\G \to S$, and the choice of a $\Gamma$
invariant section of $\U_{x_o,\omit} \to X'$. We will now explain
how Hodge theory provides canonical choices of both.

It follows from (\ref{hodge-lie}) that $\Gr^W_0\g \cong \s$.
Consequently, there is a canonical splitting of the homomorphism
$\g \to \s$. This induces a canonical
splitting of the homomorphism $\G \to S$, and therefore a
canonical action of $\Gamma$ on $\U_{x_o,\omit}$ and a canonical
identification $\G \cong S \ltimes \U$.

It remains to explain why there is a $\Gamma$ equivariant
section of $\U_{x_o,\omit}$. This is an elaboration of the
argument in the unipotent case.

For each $b\in X$, Hodge theory provides canonical ring isomorphisms
$$
\O(G_{x_o,b}) \cong \bigoplus_{m\ge 0} Gr^W_m \O(\G_{x_o,b})
$$
and
$$
\O(\cP_{x_o,b}) \cong \Gr^W_0 \O(\G_{x_o,b}).
$$
Moreover, it follows from (\ref{pre_var}) that these identifications
depend smoothly on $b$. Consequently, there is a smooth section $\sigma$
of the canonical projection $\G_{x_o,\omit} \to \cP_{x_o,\omit}$.
Restricting to $X'$, we obtain a canonical smooth section of the
projection $\U_{x_o,\omit} \to X'$.

\begin{proposition}
This section is $\Gamma$ equivariant.
\end{proposition}

\begin{proof}
For each $x\in X$ we have the action $\G\times \G_{x_o,x}\to\G_{x_o,x}$.
By (\ref{gpoid_dr}) the corresponding map of coordinate rings is a
morphism of MHS. By the choice of splitting of $\G \to S$, the action
of $S$ given by the splitting preserves the canonical isomorphism
$$
\O(\G_{x_o,x}) \cong \bigoplus_{l\ge 0}\Gr^W_l \O(\G_{x_o,x}).
$$
It follows that the section of $\G_{x_o,\omit} \to \cP_{x_o,\omit}$
defined above is equivariant with respect to the left $S$ actions.
It follows that the restriction of this section to $X'$ is $\Gamma$
equivariant.
\end{proof}

\begin{example}
In this example, we take $X$ to be the complement in $\C^n$ of
the universal discriminant locus. (View $\C^n$ as the space of
monic polynomials of degree $n$.) Pick a base point $x_o$. The
fundamental group of this
space is the classical braid group. Denote the symmetric
group on $n$ letters by $\Sigma_n$. There is a natural
homomorphism $\rho : B_n \to \Sigma_n$. Denote the
corresponding covering of $X$ by $\pi : X' \to X$. Its fundamental
group is the pure braid group $P_n$. As is
well known, $X'$ is the complement of the hyperplanes
$x_i = x_j$ in $\C^n$ where $i\neq j$. The projection takes
$(x_1,\dots,x_n)$ to the monic polynomial $\prod(T-x_j)$. The
natural left action of $\Sigma_n$ on $X'$ is given by
$$
\sigma : (x_1,\dots,x_n) \mapsto
(x_{\sigma^{-1}(1)},\dots,x_{\sigma^{-1}(n)}).
$$

The local system $\pi_\ast\Q_{X'}$ is an admissible variation of
Hodge structure over $X$ of weight 0, rank $n$, and type $(0,0)$.
The closure of the image of the monodromy is $\Sigma_n$, a
semi-simple group. So we can
apply Theorem~\ref{extended} to deduce the existence of a MHS
on the relative completion, and the existence of a universal
connection. In this case, the canonical connection is well known
by the work \cite{kohno} of Kohno.

Denote the free Lie algebra over $\C$ generated by the $Y_j$ by
$\bL(Y_1,\dots,Y_m)$. Denote the unipotent completion of $P_n$
by $\P_n$ and its Lie algebra by $\p_n$.
The associated graded of $\p_n$ of is the graded Lie algebra
$$
\bL(X_{ij} : ij\text{ is a two element subset of }\{1,\dots,n\})/R
$$
where $R$ is the ideal generated by the quadratic relations
\begin{align*}\label{braid_relns}
[X_{ij},X_{kl}]&\text{ when $i,j,k$ and $l$ are distinct;}\cr
[X_{ij},X_{ik} + X_{jk}]& \text{ when $i,j$ and $k$ are distinct}.
\end{align*}
The natural (left) action of the symmetric group on it is defined by
$$
Ad(\sigma): X_{ij} \mapsto X_{\sigma(ij)}.
$$

The canonical invariant form
$$
\omega \in E^1(X')\otimes \Gr^W_\dt \p_n
$$
is
$$
\omega = \sum_{ij} d\log(x_i - x_j) X_{ij}.
$$
It is invariant because
$$
\sigma^\ast \omega =
\sum_{ij} d\log(x_{\sigma^{-1}(i)} - x_{\sigma^{-1}(j)}) X_{ij}
= \sum_{ij} d\log(x_i - x_j) X_{\sigma(ij)}
= Ad(\sigma)\omega.
$$
We therefore obtain a homomorphism
$$
B_n \to \Sigma_n\ltimes \cP_n
$$
where $\cP_n$ denotes the complex form of the Malcev completion of
$P_n$. This is the completion of $B_n$ relative to
$\rho : B_n \to \Sigma_n$.
\end{example}

\end{document}